# THE SPHERICAL CONVEX FLOATING BODY

FLORIAN BESAU & ELISABETH M. WERNER

**Abstract**

For a convex body on the Euclidean unit sphere the spherical convex floating body is introduced. The asymptotic behavior of the volume difference of a spherical convex body and its spherical floating body is investigated. This gives rise to a new spherical area measure, the floating area. Remarkably, this floating area turns out to be a spherical analogue to the classical affine surface area from affine differential geometry. Several properties of the floating area are established.

## 1. Introduction

The floating body appeared first in the work of C. Dupin [15] in 1822. By the end of the 20th century this basic notion witnessed a surge in interest. In 1990, a seminal new definition was given by C. Schütt and E. Werner [71]. They introduced the *convex* floating body as the intersection of all halfspaces whose hyperplanes cut off a set of fixed volume of a convex body (compact convex set). In contrast to the original definition, the convex floating body is always convex and coincides with Dupin's floating body if it exists.

The main reason that the (convex) floating body has attracted considerable interest in recent decades is that it allows extensions of the classical notion of affine surface area to general convex bodies in all dimensions. Indeed, as was shown by K. Leichtweiß [36] and C. Schütt and E. Werner [71], the affine surface area arises as a limit of the volume difference of the convex body and its floating body.

Affine surface area was introduced by W. Blaschke [10] in 1923 for smooth convex bodies in Euclidean space of dimensions 2 and 3 and for smooth convex bodies. Even though it proved much more difficult to extend affine surface area to general convex bodies than other notions, like surface area measures or curvature measures, successively such extensions were achieved. Aside from the afore mentioned successful approach via the (convex) floating body, E. Lutwak [44] was able to provide an extension in 1991 by a completely different method and also proved the long conjectured upper semicontinuity of affine surface area.

As the name suggests, affine surface area is invariant under volume preserving affine transformations. Furthermore it is a valuation on the space of convex bodies and, as mentioned, upper semicontinuous. M. Ludwig and M. Reitzner [41] proved that these three properties essentially characterize affine surface area. They showed that a valuation on convex bodies that is upper semicontinuous and invariant under volume preserving affine transformations is a linear combination

*Mathematics Subject Classification.* 52A55 (Primary), 28A75, 52A20, 53A35 (Secondary).
*Key words and phrases.* spherical convexity, spherical convex floating body, floating area.
First author is supported by the European Research Council; Project number: 306445.
Second author is partially supported by an NSF grant.

of affine surface area, volume, and the Euler characteristic. Building on results of M. Ludwig and M. Reitzner [**42**], this result was recently considerably strengthened by C. Haberl and L. Parapatits [**29**].

Affine surface area is among the most powerful tools in equiaffine differential geometry (see B. Andrews [**5**, **6**], A. Stancu [**75**, **76**], M. Ivaki [**32**] and M. Ivaki and A. Stancu [**33**]). It appears naturally as the Riemannian volume of a smooth convex hypersurface with respect to the affine metric (or Berwald-Blaschke metric), see e.g. the thorough monograph of K. Leichtweiß [**37**] or the book by K. Nomizu and T. Sasaki [**55**]. In particular the upper semicontinuity proved to be critical in the solution of the affine Plateau problem by N. S. Trudinger and X.-J. Wang [**77**].

A variant of the convex floating body provided a geometric interpretation of $L_p$-affine surface area, see C. Schütt and E. Werner [**73**]. $L_p$-affine surface area is a generalization of affine surface area in the $L_p$-Brunn-Minkowski theory introduced by E. Lutwak [**45**] (also see D. Hug [**31**] and M. Meyer and E. Werner [**54**]). M. Ludwig and M. Reitzner [**42**] recently generalized $L_p$-affine surface area to Orlicz affine surface area.

One of the fundamental inequalities for affine surface area is the affine isoperimetric inequality (see W. Blaschke [**10**], L. A. Santaló [**60**] and C. M. Petty [**58**]) which states that, among all convex bodies with fixed volume, ellipsoids have the largest affine surface area. This inequality is related to various other inequalities, see E. Lutwak [**43**] and E. Lutwak, D. Yang and G. Zhang [**47**]. In particular, the affine isoperimetric inequality implies the Blaschke-Santaló inequality and it proved to be the key ingredient in the solution of many problems, see e.g. the book by R. J. Gardner [**22**] and also [**13**, **21**, **30**, **40**, **46**, **67**, **75**, **78**].

Applications of affine surface areas have been manifold. For instance, affine surface area appears in best and random approximation of convex bodies by polytopes, see K. Böröczky Jr. [**11**, **12**], P. M. Gruber [**26**–**28**], M. Ludwig [**38**], M. Reitzner [**59**] and C. Schütt [**68**, **70**]. See also C. Schütt and E. Werner [**72**] where they show that random approximation is almost as good as best approximation.

Furthermore, recent contributions indicate astonishing developments which open up new connections of affine surface area to, e.g., concentration of volume (e.g. [**19**, **42**, **48**]), and information theory (e.g. [**7**, **14**, **49**–**51**, **57**, **78**, **79**]).

In this paper we introduce the new notion of *spherical convex floating body* of a convex body on the Euclidean unit sphere. While there are important, intermittent results in the theory of spherical convex bodies, see e.g. [**18**, **61**, **62**], research has now picked up momentum. In particular the integral geometry of convex bodies in more general spaces, including the sphere, has seen tremendous progress over the last decades, see e.g. [**1**–**3**, **8**, **25**, **34**, **63**, **66**, **74**], and in recent years classical spherical convex geometry has made progress at a rapid pace, see [**4**, **9**, **17**, **20**, **23**, **24**, **53**, **56**, **64**, **80**].

Our main contributions are to introduce and investigate the spherical convex floating body and to establish the asymptotic behavior of the volume difference of a spherical convex body and its floating body. We prove that this volume difference gives rise to a new spherical area measure, the floating area, which bears striking similarities to the Euclidean affine surface area. Because of these similarities we are convinced that the floating area will become a powerful tool in spherical convex geometry.

The next section gives the precise statement of our main theorem. In Section 3 we will briefly recall basic notions from Euclidean convex geometry. In particular, we recall results about the convex floating body and affine surface area. In Section 4 we will fix the notation for spherical convex geometry and investigate the boundary structure of a spherical convex body. The spherical convex floating body is introduced in Section 5. We proof the main theorem in Section 6. In the final Section 7 we investigate the floating area, or more generally, the floating measure. We show that the floating area is an upper semicontinuous valuation and establish a duality formula which involves the polar body. Moreover, in Subsection 7.1. we derive an isoperimetric inequality for the floating area.

## 2. Statement of the Main Theorem

We denote the Euclidean unit sphere in $\mathbb{R}^{n+1}$ by $\mathbb{S}^n$, $n \geq 2$. The natural spherical distance $d$ is given by $d(u,v) = \arccos(u \cdot v)$ for $u, v \in \mathbb{S}^n$, where $\cdot$ denotes the Euclidean scalar product. A set $A \subseteq \mathbb{S}^n$ is called *(spherical) convex* if its radial extension

$$\operatorname{rad} A = \{rv \in \mathbb{R}^{n+1} : r \geq 0 \text{ and } v \in A\} \tag{2.1}$$

is convex in $\mathbb{R}^{n+1}$.

A closed convex subset of $\mathbb{S}^n$ is called a *(spherical) convex body*. If a convex body does not contain a pair of antipodal points, we call it *proper*. The set of convex bodies (resp. proper convex bodies) is denoted by $\mathcal{K}(\mathbb{S}^n)$ (resp. $\mathcal{K}^p(\mathbb{S}^n)$). Furthermore, the subset of bodies with non-empty interior is denoted by $\mathcal{K}_0(\mathbb{S}^n)$ (resp. $\mathcal{K}_0^p(\mathbb{S}^n)$).

A $k$-sphere $S$, $0 \leq k \leq n$, is the intersection of a $(k+1)$-dimensional linear subspace $L$ of $\mathbb{R}^{n+1}$ with $\mathbb{S}^n$. For $u \in \mathbb{S}^n$ we denote by $\mathbb{S}_u$ the hypersphere that is given by $\mathbb{S}_u = \{v \in \mathbb{S}^n : u \cdot v = 0\}$. The closed hemisphere with center in $u$ is denoted by $\mathbb{S}_u^+$. We set $\mathbb{S}_u^- = \mathbb{S}_{-u}^+$.

The $k$-dimensional Hausdorff measure is denoted by $\mathcal{H}^k$ and $\operatorname{vol}_n$ denotes the usual volume measure on $\mathbb{S}^n$. Note that $\operatorname{vol}_n$ coincides with the $n$-dimensional Hausdorff measure restricted to $\mathbb{S}^n$.

In Euclidean convex geometry the convex floating body is defined by the intersection of all halfspaces such that the hyperplanes cut off a set of constant volume. This motivates the following

**Definition** (Spherical Convex Floating Body)**.** *Let $K \in \mathcal{K}_0(\mathbb{S}^n)$, $K \neq \mathbb{S}^n$. For $\delta > 0$ small enough we define the* (spherical) convex floating body $K_{[\delta]}$ *as the intersection of all closed hemispheres $\mathbb{S}^-$ such that the hyperspheres cut off a set of volume less or equal $\delta$, that is,*

$$K_{[\delta]} = \bigcap \{\mathbb{S}^- : \operatorname{vol}_n(K \cap \mathbb{S}^+) \leq \delta\}, \tag{2.2}$$

*where $\mathbb{S}^+$ is the complementary closed hemisphere to $\mathbb{S}^-$, that is, $\mathbb{S}^+ = -\mathbb{S}^-$.*

In our main theorem we consider the volume difference of a spherical convex body $K$ and its floating body. We show that the limit, as $\delta$ goes to zero, converges to the total curvature over the boundary $\operatorname{bd} K$ of $K$ when integrating the spherical Gauss-Kronecker curvature $H_{n-1}^{\mathbb{S}^n}(K, \cdot)$ raised to the power $\frac{1}{n+1}$ (see Section 4 for details on the spherical Gauss-Kronecker curvature).

**Theorem 2.1.** *If $K \in \mathcal{K}_0^p(\mathbb{S}^n)$ is a proper convex body with non-empty interior, then*

$$\lim_{\delta \to 0^+} \frac{\mathrm{vol}_n(K) - \mathrm{vol}_n(K_{[\delta]})}{\delta^{\frac{2}{n+1}}} = c_n \int_{\mathrm{bd}\, K} H_{n-1}^{\mathbb{S}^n}(K,w)^{\frac{1}{n+1}} \, d\mathcal{H}^{n-1}(w), \qquad (2.3)$$

*where $c_n = \frac{1}{2}\left(\frac{n+1}{\kappa_{n-1}}\right)^{\frac{2}{n+1}}$ and $\kappa_{n-1}$ is the volume of the $(n-1)$-dimensional Euclidean unit ball.*

We will prove Theorem 2.1 in Section 6.

The curvature integral appearing on the right-hand side of (2.3) shares striking similarities with the affine surface area from affine differential geometry, see (3.2). Since, to our knowledge, there is no property similar to the affine invariance of Euclidean affine surface area in the spherical setting, we are reluctant to call this new spherical area measure a spherical affine surface area.

**Definition** (Floating Area). *For a convex body $K \in \mathcal{K}_0(\mathbb{S}^n)$ with non-empty interior the floating area $\Omega(K)$ is defined by*

$$\Omega(K) = \int_{\mathrm{bd}\, K} H_{n-1}^{\mathbb{S}^n}(K,\cdot)^{\frac{1}{n+1}} \, d\mathcal{H}^{n-1},$$

*if $K$ is proper or $\Omega(K) = 0$ otherwise.*

In Section 7 we investigate properties of the floating area.

## 3. The Floating Body and Affine Surface Area

In this section we collect well known result from Euclidean convex geometry. For a general reference on these facts we refer to [65]. We denote by $\cdot$ the scalar product and by $\|.\|$ the Euclidean norm. A convex body $K$ is a compact convex subset of $\mathbb{R}^n$. It is uniquely determined by its *support function* $h_K$ defined by

$$h_K(x) = \max\{x \cdot y : y \in K\}, \quad x \in \mathbb{R}^n.$$

The set of convex bodies in $\mathbb{R}^n$ is denoted by $\mathcal{K}(\mathbb{R}^n)$. We denote by $\mathcal{K}_0(\mathbb{R}^n)$ the set of convex bodies with non-empty interior.

**Definition** (Convex Floating Body [71]). *Let $K \in \mathcal{K}_0(\mathbb{R}^n)$. For $\delta > 0$ the convex floating body $K_{[\delta]}$ is defined as the intersection of all closed halfspaces $\mathbb{H}^-$ such that the hyperplanes cut off a set of volume less or equal $\delta$, that is,*

$$K_{[\delta]} = \bigcap \{\mathbb{H}^- : \mathrm{vol}_n(K \cap \mathbb{H}^+) \leq \delta\}.$$

Basic properties of the convex floating body are collected in the following

**Proposition 3.1** ([71]). *Let $K \in \mathcal{K}_0(\mathbb{R}^n)$ and $\delta > 0$ such that $K_{[\delta]}$ exists.*
  (i) *Through every $x \in \mathrm{bd}\, K_{[\delta]}$ there exists at least one hyperplane $\mathbb{H}$ that cuts off of $K$ a set of volume $\delta$.*
  (ii) *A hyperplane $\mathbb{H}$ that cuts off of $K$ a set of volume $\delta$ touches $K_{[\delta]}$ in exactly one point, the barycenter of $K \cap \mathbb{H}$.*
 (iii) *$K_{[\delta]}$ is strictly convex.*
 (iv) *Let $\delta_0 = \max\{\delta : \mathrm{vol}_n(K_{[\delta]}) > 0\}$. Then $K_{[\delta_0]}$ is only one point and for all $0 < \delta < \delta_0$, $K_{[\delta]}$ exists and has non-empty interior.*
  (v) *For a linear transformation $A \in \mathrm{GL}(\mathbb{R}^n)$ and a vector $y \in \mathbb{R}^n$ we have $(AK + y)_{[\delta]} = AK_{[|\det A|\,\delta]} + y.$*

We can parametrize the halfspaces in the definition of the convex floating body over $\mathbb{S}^{n-1}$ in the following way:

**Corollary 3.2.** *If $K \in \mathcal{K}_0(\mathbb{R}^n)$, then for any $v \in \mathbb{S}^{n-1}$ there exists $s(v, \delta) \in \mathbb{R}$ such that*
$$\delta = \operatorname{vol}_n\left(K \cap \mathbb{H}^+_{v,s(v,\delta)}\right),$$
*where $\mathbb{H}^+_{v,s(v,\delta)} = \{x \in \mathbb{R}^n : x \cdot v \geq s(v,\delta)\}$. Moreover, we have*
$$K_{[\delta]} = \bigcap_{v \in \mathbb{S}^{n-1}} \mathbb{H}^-_{v,s(v,\delta)}. \tag{3.1}$$

*Proof.* This follows from Proposition 3.1 (i). q.e.d.

The generalized Gauss-Kronecker curvature of a convex body $K \in \mathcal{K}(\mathbb{R}^n)$ at a boundary point $x$ is denoted by $H^{\mathbb{R}^n}_{n-1}(K, x)$ and exists for $\mathcal{H}^{n-1}$-almost all boundary points, see e.g. [**31**, Lemma 2.3].

**Definition** (Affine Surface Area). *Let $K \in \mathcal{K}_0(\mathbb{R}^n)$. Then the affine surface area $\operatorname{as}(K)$ of $K$ is defined by*
$$\operatorname{as}(K) = \int_{\operatorname{bd} K} H^{\mathbb{R}^n}_{n-1}(K, x)^{\frac{1}{n+1}} dx. \tag{3.2}$$

Affine surface area $\operatorname{as}(K)$ is finite for all convex bodies $K$. This can be seen in the following way: For a convex body $K \in \mathcal{K}_0(\mathbb{R}^n)$ and a boundary point $x \in \operatorname{bd} K$ we denote by $r_K(x)$ the maximal radius of a Euclidean ball that is contained in $K$ and touches the boundary of $K$ in $x$, in other words,
$$r_K(x) = \sup\{r \geq 0 : \exists y \in K \text{ such that } B^n(y, r) \subseteq K \text{ and } x \in \operatorname{bd} B^n(y, r)\},$$
where $B^n(y, r)$ denotes the closed Euclidean ball of radius $r$ with center $y$.

By Blaschke's Rolling Theorem (see e.g. [**65**, Corollary 3.2.13]) we know that $r_K > 0$ for $\mathcal{H}^{n-1}$-almost all boundary points. In fact C. Schütt and E. Werner proved in [**71**] the following:

**Theorem 3.3** ([**71**]). *Let $K \in \mathcal{K}_0(\mathbb{R}^n)$. Then for all $\alpha \in [0, 1)$*
$$\int_{\operatorname{bd} K} r_K(x)^{-\alpha} dx < \infty. \tag{3.3}$$

Since the Gauss-Kronecker curvature at a boundary point $x$ of $K$ is less or equal to the curvature of any ball contained in $K$ that touches the boundary in $x$ we have $H^{\mathbb{R}^n}_{n-1}(K, x) \leq r_K(x)^{-(n-1)}$ for $\mathcal{H}^{n-1}$-almost all $x \in \operatorname{bd} K$. We conclude
$$\operatorname{as}(K) \leq \int_{\operatorname{bd} K} r_K(x)^{-\frac{n-1}{n+1}} dx < \infty.$$

The limit of the volume difference of a convex body and its floating body converges to the affine surface area of the body in the following way:

**Theorem 3.4** ([**71**]). *Let $K \in \mathcal{K}_0(\mathbb{R}^n)$. Then*
$$\operatorname{as}(K) = \frac{1}{c_n} \lim_{\delta \to 0^+} \frac{\operatorname{vol}_n(K) - \operatorname{vol}_n(K_{[\delta]})}{\delta^{\frac{2}{n+1}}}, \tag{3.4}$$
*where $c_n = \frac{1}{2}\left(\frac{n+1}{\kappa_{n-1}}\right)^{\frac{2}{n+1}}$.*

By (3.4) and the covariance of the floating body under affine transformations that preserve volume we conclude that for all $A \in \mathrm{SL}(\mathbb{R}^n)$ and $y \in \mathbb{R}^n$

$$\mathrm{as}(AK + y) = \mathrm{as}(K).$$

We see that, as the name suggests, affine surface area is invariant under volume preserving affine transformations.

The proof of Theorem 3.4 is built on the following results. We choose to include them here since we will need them to prove Theorem 2.1. Note that we denote the convex hull of two points $x, y \in \mathbb{R}^n$ by $\mathrm{conv}(x, y)$.

**Theorem 3.5** ([71]). *Let $K \in \mathcal{K}_0(\mathbb{R}^n)$ and $0 \in \mathrm{int}\, K$. For $x \in \mathrm{bd}\, K$ and $\delta > 0$ we set $\{x_\delta\} = \mathrm{bd}\, K_{[\delta]} \cap \mathrm{conv}(0, x)$. Then there exists $C > 0$ and $\delta_0 > 0$ such that for all $\delta < \delta_0$ we have*

$$\frac{\|x - x_\delta\|}{\delta^{\frac{2}{n+1}}} \leq C r_K(x)^{-\frac{n-1}{n+1}},$$

*for $\mathcal{H}^{n-1}$-almost all $x \in \mathrm{bd}\, K$.*

**Theorem 3.6** ([71]). *Let $K \in \mathcal{K}_0(\mathbb{R}^n)$ with $0 \in \mathrm{int}\, K$. For $\delta > 0$ small enough we have for $\mathcal{H}^{n-1}$-almost all $x \in \mathrm{bd}\, K$*

$$\lim_{\delta \to 0^+} \frac{1}{n} \frac{x \cdot N_x^K}{\delta^{\frac{2}{n+1}}} \left(1 - \left(\frac{\|x_\delta\|}{\|x\|}\right)^n\right) = c_n H_{n-1}^{\mathbb{R}^n}(K, x)^{\frac{1}{n+1}}, \qquad (3.5)$$

*where $\{x_\delta\} = \mathrm{bd}\, K_{[\delta]} \cap \mathrm{conv}(0, x)$ and $N_x^K$ denotes the outer unit normal vector at $x \in \mathrm{bd}\, K$.*

Note that for the left hand side in (3.5) we have

$$\lim_{\delta \to 0^+} \frac{1}{n} \frac{x \cdot N_x^K}{\delta^{\frac{2}{n+1}}} \left(1 - \left(\frac{\|x_\delta\|}{\|x\|}\right)^n\right) = \lim_{\delta \to 0^+} \left(\frac{x}{\|x\|} \cdot N_x^K\right) \frac{\|x - x_\delta\|}{\delta^{\frac{2}{n+1}}}. \qquad (3.6)$$

## 4. Basic Facts from Spherical Convex Geometry

The *convex hull* $\mathrm{conv}\, A$ of $A \subseteq \mathbb{S}^n$ is the intersection of all convex bodies in $\mathbb{S}^n$ that contain $A$, with the usual convention that the empty intersection is the whole sphere $\mathbb{S}^n$. The convex hull of two convex bodies $K, L \in \mathcal{K}(\mathbb{S}^n)$ is denoted by $\mathrm{conv}(K, L) = \mathrm{conv}(K \cup L)$ and the segment spanned by two points $u, v \in \mathbb{S}^n$, $u \neq -v$, is given by $\mathrm{conv}(u, v) = \mathrm{conv}(\{u\}, \{v\})$.

The interior of a set $A$ is denoted by $\mathrm{int}\, A$. We denote a closed spherical cap of radius $0 \leq \alpha \leq \frac{\pi}{2}$ and center $u \in \mathbb{S}^n$ by $C_u(\alpha) = \{v \in \mathbb{S}^n : d(u, v) \leq \alpha\}$.

The *polar body* $K^\circ$ of a convex body $K \in \mathcal{K}(\mathbb{S}^n)$ is defined by

$$K^\circ = \{v \in \mathbb{S}^n : v \cdot w \leq 0 \text{ for all } w \in K\},$$

and is again a convex body. The following lemma collects well-known facts about the polar body:

**Lemma 4.1.** *Let $K \in \mathcal{K}(\mathbb{S}^n)$. Then*
*(i) $(K^\circ)^\circ = K$.*
*(ii) We have*

$$K^\circ = \bigcap_{w \in K} \mathbb{S}_w^- = \{v : K \subseteq \mathbb{S}_v^-\}.$$

*In particular,*
$$\text{int } K^\circ = \bigcap_{w \in K} \text{int } \mathbb{S}_w^- = \{v : K \subseteq \text{int } \mathbb{S}_v^-\} = \mathbb{S}^n \setminus K_{\frac{\pi}{2}}$$

*where $K_{\frac{\pi}{2}}$ is the set of all points that have distance of at most $\frac{\pi}{2}$ to $K$.*

(iii) *If $K$ is proper, then $\text{int } K = \emptyset$ if and only if $\text{int } K^\circ = \emptyset$.*

(iv) $(K \cap L)^\circ = \text{conv}(K^\circ, L^\circ).$

(v) *If $K \in \mathcal{K}_0^p(\mathbb{S}^n)$, then $\text{int } K \cap \text{int}(-K^\circ) \neq \emptyset$. In particular, we have for all $u \in \text{int } K \cap \text{int}(-K^\circ)$ that $u \in \text{int } K$ and $K \subseteq \text{int } \mathbb{S}_u^+$.*

Let $S$ be a $k$-sphere for some $k \in \{0, \ldots, n\}$ and let $K \in \mathcal{K}(\mathbb{S}^n)$. Then the *spherical projection* $K|S$ is defined by
$$K|S = \text{conv}(K, S^\circ) \cap S.$$

Note that, since $S$ is a $k$-sphere, $S^\circ$ is an $(n-k-1)$-sphere.

The spherical projection of a point $w$ to a hypersphere $\mathbb{S}_v$, for $v \neq \pm w$, is given by $w|\mathbb{S}_v = \frac{w - \cos(d(v,w))v}{\sin(d(v,w))}$. We have $w = \cos(d(v,w))v + \sin(d(v,w))(w|\mathbb{S}_v)$.

**4.1. Spherical Support Function and Gnomonic Projection.** For the results stated here see, e.g., [9].

**Definition** (Spherical Support Function). *For a fixed $u \in \mathbb{S}^n$ and a proper convex body $K \in \mathcal{K}^p(\mathbb{S}^n)$ such that $K \subseteq \text{int } \mathbb{S}_u^+$, the spherical support function $h_u(K, \cdot) \colon \mathbb{S}_u \to \left(-\frac{\pi}{2}, \frac{\pi}{2}\right)$ of $K$ is defined by*
$$h_u(K, v) = \max\{\text{sgn}(v \cdot w) d(u, w|\mathbb{S}_{u,v}^1) : w \in K\},$$
*where $\mathbb{S}_{u,v}^1$ denotes the 1-sphere spanned by $u$ and $v$.*

The intuitive interpretation of the spherical support function is as follows: If $v \in \mathbb{S}_u$ then the projection $K|\mathbb{S}_{u,v}^1$ is a spherical segment and the spherical support function measures the width with respect to $u$. We have
$$K|\mathbb{S}_{u,v}^1 = \{\cos(\alpha)u + \sin(\alpha)v : \alpha \in [-h_u(K, -v), h_u(K, v)]\}.$$

Let $u \in \mathbb{S}^n$. For $v \in \mathbb{S}_u$ and $\delta \in \left(-\frac{\pi}{2}, \frac{\pi}{2}\right)$, we set $z = \cos(\delta)v - \sin(\delta)u$ and
$$\mathbb{S}_{u,v,\delta}^+ = \mathbb{S}_z^+ = \{w \in \mathbb{S}^n : v \cdot w \geq \tan(\delta)(u \cdot w)\}. \tag{4.1}$$

Let $K \in \mathcal{K}(\mathbb{S}^n)$ such that $K \subseteq \text{int } \mathbb{S}_u^+$. For $v \in \mathbb{S}_u$, the the supporting hyperspheres of $K$ are parametrized by $\mathbb{S}_{u,v,h_u(K,v)}$. In particular, we have that $z = \cos(h_u(K,v))v - \sin(h_u(K,v))u \in \text{bd } K^\circ$ and $z$ is an outer unit normal vector to $K$ for every point $w \in \text{bd } K \cap \mathbb{S}_{u,v,h_u(K,v)}$. We note that $v = z|\mathbb{S}_u$ and
$$-u \cdot z = \sin(h_u(K, z|\mathbb{S}_u)). \tag{4.2}$$

We conclude that $K = \bigcap_{v \in \mathbb{S}_u} \mathbb{S}_{u,v,h_u(K,v)}^- = \bigcap_{z \in \text{bd } K^\circ} \mathbb{S}_z^-$.

**Definition** (Gnomonic Projection). *For $u \in \mathbb{S}^n$ we denote by $\mathbb{R}_u^n$ the linear subspace orthogonal to $u$ in $\mathbb{R}^{n+1}$, i.e., $\mathbb{R}_u^n = \{x \in \mathbb{R}^{n+1} : u \cdot x = 0\}$. The gnomonic projection $g_u \colon \text{int } \mathbb{S}_u^+ \to \mathbb{R}_u^n$ is defined by*
$$g_u(v) = \frac{v}{u \cdot v} - u.$$

The inverse of the gnomonic projection $g_u^{-1} \colon \mathbb{R}_u^n \to \text{int } \mathbb{S}_u^+$ is
$$g_u^{-1}(x) = \frac{x + u}{\|x + u\|}, \quad x \in \mathbb{R}_u^n.$$

The geometric interpretation of the gnomonic projection is as follows: If $A \subseteq \mathbb{S}_u^+$ then the image of $A$ under the gnomonic projection is given by the radial extension of $A$ and the intersection with the tangent plane $u + T_u \mathbb{S}^n$ (for convenience we project $u + T_u \mathbb{S}^n$ to $\mathbb{R}_u^n$ in the end). Thus, if $S$ is a $k$-sphere for some $k \in \{0, \ldots, n-1\}$, then $g_u(S \cap \operatorname{int} \mathbb{S}_u^+)$ is an affine $k$-dimensional subspace of $\mathbb{R}_u^n$. In particular, $g_u$ maps convex bodies in $\operatorname{int} \mathbb{S}_u^+$ bijectively to $\mathcal{K}(\mathbb{R}_u^n)$. For the spherical support function of $K$ we have (see, e.g. [9])

$$\tan h_u(K, v) = h_{g_u(K)}(v), \quad v \in \mathbb{S}_u, \tag{4.3}$$

where $h_{g_u(K)} \colon \mathbb{R}_u^n \to \mathbb{R}$ denotes the Euclidean support function of $g_u(K)$ in $\mathbb{R}_u^n$.

For $K \in \mathcal{K}_0^p(\mathbb{S}^n)$ and $u \in \mathbb{S}^n$ such that $u \in \operatorname{int} K$, the *radial function* $\rho_u(K, \cdot) \colon \mathbb{S}_u \to [0, \pi]$ is defined by

$$\rho_u(K, v) = \max\{\alpha \in [0, \pi] : \cos(\alpha) u + \sin(\alpha) v \in K\}. \tag{4.4}$$

If $K \subseteq \operatorname{int} \mathbb{S}_u^+$, then the Euclidean radial function $\rho_{g_u(K)} \colon \mathbb{R}_u^n \to \mathbb{R}$ is given by

$$\tan(\rho(K, v)) = \rho_{g_u(K)}(v), \quad v \in \mathbb{S}_u.$$

The Euclidean polar body $g_u(K)^* = \{y \in \mathbb{R}_u^n : x \cdot y \leq 1 \text{ for all } x \in g_u(K)\}$ is related to $K^\circ$ by $g_u(K)^* = g_{-u}(K^\circ)$. Moreover,

$$\tan(h_u(K, v)) = h_{g_u(K)}(v) = \frac{1}{\rho_{g_u(K)^*}(v)} = \cot(\rho_{-u}(K^\circ, v)),$$

or, equivalently,

$$h_u(K, v) + \rho_{-u}(K^\circ, v) = \frac{\pi}{2}, \quad v \in \mathbb{S}_u. \tag{4.5}$$

**4.2. Boundary Structure of Spherical Convex Bodies.** In this section we develop the technical framework to transform integrals on $\mathbb{S}^n$ and the boundary of spherical convex bodies via the gnomonic projection or the Gauss map. We will consider subsets of the sphere as subsets in $\mathbb{R}^{n+1}$ and use the area formula on rectifiable sets, where we explicitly calculate the (approximate) tangential Jacobian. For a reference on the area formula and tangential Jacobian we refer to F. Maggi [52] or S. G. Krantz and H. R. Parks [35], which will provide sufficient background for the tools we use (for a more extensive position see, e.g., H. Federer [16]).

We begin with an outline of what follows: Using the area formula (see Theorem 11.6 in [52]) we show for the diffeomorphism $g_u \colon \operatorname{int} \mathbb{S}_u^+ \to \mathbb{R}_u^n$, that for measurable $\omega \subseteq \operatorname{int} \mathbb{S}_u^+$ and measurable $f \colon \mathbb{R}_u^n \to \mathbb{R}$ we have

$$\int_\omega f \circ g_u \, d\mathcal{H}_{\mathbb{S}^n}^{n-1} = \int_{g_u(\omega)} f J^{\operatorname{int} \mathbb{S}_u^+}(g_u) \, d\mathcal{H}_{\mathbb{R}_u^n}^{n-1}.$$

Here $J^{\operatorname{int} \mathbb{S}_u^+}(g_u)$ denotes the *tangential Jacobian* of $g_u$ on $\operatorname{int} \mathbb{S}_u^+$ which is defined by

$$J^{\operatorname{int} \mathbb{S}_u^+}(g_u)(v) = \sqrt{\det((d(g_u)_v)^* d(g_u)_v)},$$

where $d(g_u)_v$ denotes the *tangential derivative* (differential) in $v$ and $(d(g_u)_v)^*$ denotes the adjoint. In the following proposition we will explicitly calculate this expression.

**Proposition 4.2.** *Let $u \in \mathbb{S}^n$. Then $g_u$ is a diffeomorphism between $\operatorname{int} \mathbb{S}_u^+$ and $\mathbb{R}_u^n$. For $v \in \operatorname{int} \mathbb{S}_u^+$ we identify the tangent space $T_v \operatorname{int} \mathbb{S}_u^+$ with $\mathbb{R}_v^n$ and* $T_{g_u(v)}\mathbb{R}^n_u$ with $\mathbb{R}^n_u$. Then the differential $d(g_u)_v \colon T_v\mathrm{int}\,\mathbb{S}^+_u \to T_{g_u(v)}\mathbb{R}^n_u$ is given by

$$d(g_u)_v(X_v) = \frac{1}{u \cdot v} X_v - \frac{u \cdot X_v}{(u \cdot v)^2} v, \tag{4.6}$$

for all $X_v \in \mathbb{R}^n_v$. The tangential Jacobian of $g_u$ is given by

$$J^{\mathrm{int}\,\mathbb{S}^+_u}(g_u)(v) = \sqrt{\det((d(g_u)_v)^* d(g_u)_v)} = \frac{1}{(u \cdot v)^{n+1}} \tag{4.7}$$

and the tangential Jacobian of the inverse $g_u^{-1}$ in $x \in \mathbb{R}^n_u$ is

$$J^{\mathbb{R}^n_u}(g_u^{-1})(x) = \frac{1}{J^{\mathrm{int}\,\mathbb{S}^+_u}(g_u)(g_u^{-1}(x))} = \frac{1}{(1+\|x\|^2)^{\frac{n+1}{2}}}. \tag{4.8}$$

*Proof.* An elementary calculation leads to formula (4.6). Also the inverse can be explicitly calculated. Thus $g_u$ is a diffeomorphism. In order to prove (4.7), we first note that $d(g_u)_v$ can be expressed as a matrix on $\mathbb{R}^{n+1}$ of rank $n$ by

$$d(g_u)_v = \frac{1}{u \cdot v}\left(\mathrm{Id}_{n+1} - \frac{1}{u \cdot v} v \otimes u\right),$$

where $v \otimes u$ denotes the matrix determined by $(v \otimes u)X_v = (u \cdot X_v)v$. Thus,

$$(d(g_u)_v)^* d(g_u)_v = \frac{1}{(u \cdot v)^2}\left(\mathrm{Id}_{\mathbb{R}^n_v} + \frac{u - (u \cdot v)v}{u \cdot v} \otimes \frac{u - (u \cdot v)v}{u \cdot v}\right).$$

Using the well known formula $\det(\mathrm{Id} + z \otimes z) = 1 + \|z\|^2$, we conclude

$$J^{\mathrm{int}\,\mathbb{S}^+_u}(g_u)(v) = \frac{1}{(u \cdot v)^n}\sqrt{1 + \frac{\|u - (u \cdot v)v\|^2}{(u \cdot v)^2}} = \frac{1}{(u \cdot v)^{n+1}}.$$

q.e.d.

The analytic properties of the boundary of a spherical convex body $K$ are similar to the properties of the boundary of a Euclidean convex body. This is obvious when considering the gnomonic projection around a boundary point $w \in \mathrm{bd}\,K$, that is, for an arbitrary but fixed $\varepsilon \in (0, \frac{\pi}{2})$ we consider $L = g_w(K \cap C_w(\varepsilon))$. Then $L$ is a Euclidean convex body with $g_w(w) = 0$ and $\mathrm{bd}\,L = \mathrm{bd}\,g_w(K \cap C_w(\varepsilon))$. Since $g_w$ is a diffeomorphism on $C_w(\varepsilon)$ we conclude that all regularity results of $\mathrm{bd}\,L$ in $0$ from Euclidean convex geometry hold for $w \in \mathrm{bd}\,K$.

**Proposition 4.3.** *Let $K \in \mathcal{K}(\mathbb{S}^n)$.*
*(i) The boundary of $K$ is an $\mathcal{H}^{n-1}$-rectifiable set.*
*(ii) If $K$ has non-empty interior, then for $\mathcal{H}^{n-1}$-almost all boundary points $w$ there exists a unique outer unit normal $N^K_w \in \mathbb{S}_w$, where we identify the tangent space $T_w\mathbb{S}^n$ with $\mathbb{R}^n_w$.*

Note that the gnomonic projection is not conformal and therefore $d(g_u)_w(N^K_w)$ is in general not the outer normal vector $N^{g_u(K)}_{g_u(w)}$ of $\mathrm{bd}\,g_u(K)$. The relation between the outer normal vector of $\mathrm{bd}\,K$ and $\mathrm{bd}\,g_u(K)$ is given by

$$N^{g_u(K)}_{g_u(w)} = N^K_w|_{\mathbb{S}_u}. \tag{4.9}$$

Now let $u \in \mathbb{S}^n$ and $K \in \mathcal{K}^p_0(\mathbb{S}^n)$ such that $K \subseteq \mathrm{int}\,\mathbb{S}^+_u$. Then there is $\beta \in (0, \frac{\pi}{2})$ such that $K \subseteq \mathrm{int}\,C_u(\beta)$. Since $g_u$ is differentiable and therefore Lipschitz on $C_u(\beta)$ and because $\mathrm{bd}\,K$ is $\mathcal{H}^{n-1}$-rectifiable, we conclude that the approximate tangential derivative $d^{\mathrm{bd}\,K}g_u$ exists $\mathcal{H}^{n-1}$-almost everywhere.

Furthermore, since $g_u$ maps tangential hypersphere to $w$ to the affine hyperplane tangent to $g_u(w)$, we have $d(g_u)_w(T_w \operatorname{bd} K) = T_{g_u(w)} \operatorname{bd} g_u(K)$ and conclude

$$d^{\operatorname{bd} K}(g_u)_w(X_w) = d(g_u)_w(X_w),$$

for all $X_w \in T_w \operatorname{bd} K$. We can write $d^{\operatorname{bd} K}(g_u)_w$ as a matrix on $\mathbb{R}^{n+1}$ of rank $n-1$ in the form

$$d^{\operatorname{bd} K}(g_u)_w = \frac{1}{u \cdot w}\left(\operatorname{Id}_{n+1} - \frac{1}{u \cdot w} w \otimes u\right) - d(g_u)_w(N_w^K) \otimes N_w^K.$$

Thus,

$$(d^{\operatorname{bd} K}(g_u)_w)^* d^{\operatorname{bd} K}(g_u)_w = \frac{1}{(u \cdot w)^2}(\operatorname{Id}_{T_w \operatorname{bd} K} + \tilde{z} \otimes \tilde{z}),$$

where $\tilde{z} = \frac{u - (u \cdot w)w - (u \cdot N_w^K)N_w^K}{u \cdot w}$. Hence, the approximate tangential Jacobian

$$J^{\operatorname{bd} K}(g_u)(w) = \sqrt{\det\left((d^{\operatorname{bd} K}(g_u)_w)^* d^{\operatorname{bd} K}(g_u)_w\right)}$$

is given by

$$J^{\operatorname{bd} K}(g_u)(w) = \frac{\sqrt{1 - (u \cdot N_w^K)^2}}{(u \cdot w)^n} = \frac{\cos(h_u(K, N_w^K|\mathbb{S}_u))}{\cos(d(u,w))^n}, \quad (4.10)$$

for $\mathcal{H}^{n-1}$-almost all $w \in \operatorname{bd} K$, where the last equation have used (4.2) for $z = N_w^K|\mathbb{S}_u$. For the inverse $g_u^{-1}$ we obtain

$$J^{\operatorname{bd} g_u(K)}(g_u^{-1})(x) = \frac{1}{J^{\operatorname{bd} K}(g_u)(g_u^{-1}(x))} = \frac{\sqrt{1 + \left(x \cdot N_x^{g_u(K)}\right)^2}}{(1 + \|x\|^2)^{\frac{n}{2}}}, \quad (4.11)$$

for $\mathcal{H}^{n-1}$-almost all $x \in \operatorname{bd} g_u(K)$.

The *radial map* $R_K \colon \operatorname{bd} K \to \mathbb{S}_u$ is defined by $R_K(w) = w|\mathbb{S}_u$. The Euclidean radial map $R_L \colon \operatorname{bd} L \to \mathbb{S}^{n-1}$ is defined by $R_L(x) = \frac{x}{\|x\|}$. The tangential Jacobian of $R_L$ is given by $J^{\operatorname{bd} L}(R_L)(x) = \frac{x \cdot N_x^L}{\|x\|^n} = \frac{h_L(N_x^L)}{\|x\|^n}$ for $\mathcal{H}^{n-1}$-almost all $x \in \operatorname{bd} L$.

The radial maps $R_K$ and $R_{g_u(K)}$ are related by $R_{g_u(K)}(g_u(w)) = R_K(w)$. We conclude for the approximate tangential Jacobian of $R_K$,

$$J^{\operatorname{bd} K}(R_K)(w) = J^{\operatorname{bd} g_u(K)}(R_{g_u(K)})(g_u(w)) J^{\operatorname{bd} K}(g_u)(w)$$
$$= \frac{\sin(h_u(K, N_w^K|\mathbb{S}_u))}{\sin(d(u,w))^n} = \frac{-u \cdot N_w^K}{\sin(d(u,w))^n}, \quad (4.12)$$

where we used the fact that $h_{g_u(K)}\left(N_{g_u(w)}^{g_u(K)}\right) = \tan(h_u(K, N_w^K|\mathbb{S}_u))$ (see (4.3), (4.2) and (4.9)).

A boundary point $w \in \operatorname{bd} K$ is called *regular* if it has a unique outer unit normal $N_w^K$ at $w$. The set of regular boundary points of $K$ is denoted by $\operatorname{reg} K$. As in the Euclidean setting, the Gauss map $N^K \colon \operatorname{reg} K \to \mathbb{S}^n$ is in general not Lipschitz, see e.g. [**31**]. However, as was pointed out in [**31**], if, for $\alpha > 0$, one restricts $N^K$ to $(\operatorname{bd} K)_\alpha$, defined by

$$(\operatorname{bd} K)_\alpha = \{w \in \operatorname{bd} K : \exists v \in \mathbb{S}^n \text{ such that } w \in C_v(\alpha) \subseteq K\},$$

then $N^K|_{(\operatorname{bd} K)_\alpha}$ is Lipschitz.

Furthermore, for
$$(\operatorname{bd} K)_+ := \bigcup_{n \in \mathbb{N}} (\operatorname{bd} K)_{\frac{1}{n}}$$

we have $\mathcal{H}^{n-1}(\operatorname{bd} K \setminus (\operatorname{bd} K)_+) = 0$, see [**31**, Lemma 2.2]. It follows that for $\mathcal{H}^{n-1}$-almost all boundary points the approximate Jacobian of $N^K$ exists and is therefore given by $J^{\operatorname{bd} K}(N^K)(w) = H_{n-1}^{\mathbb{S}^n}(K, w)$, where $H_{n-1}^{\mathbb{S}^n}(K, w)$ denotes the Gauss-Kronecker curvature. Again this can be seen easily when considering the gnomonic projection around a boundary point $w \in \operatorname{bd} K$, since

$$J^{\operatorname{bd} K}(N^K)(w) = J^{\operatorname{bd} g_w(K \cap C_w(\varepsilon))}(N^{g_w(K \cap C_w(\varepsilon))})(0) J^{\operatorname{bd} K}(g_w)(w)$$
$$= H_{n-1}^{\mathbb{R}^n_w}(g_w(K \cap C_w(\varepsilon)), 0)$$

and the fact that $H_{n-1}^{\mathbb{R}^n_w}(g_w(K \cap C_w(\varepsilon)), 0) = H_{n-1}^{\mathbb{S}^n}(K, w)$. The latter is obvious since $d(g_w)_w = \operatorname{Id}_{\mathbb{R}^n_w}$.

In particular, for a proper convex body with non-empty interior there exists $u \in \operatorname{int} K$ such that $K \subseteq \operatorname{int} \mathbb{S}_u^+$ (see Lemma 4.1(v)) and we may express $H_{n-1}^{\mathbb{S}^n}(K, \cdot)$ by $H_{n-1}^{\mathbb{R}^n_u}(g_u(K), \cdot)$.

**Lemma 4.4.** *Let $K \in \mathcal{K}_0^p(\mathbb{S}^n)$ and $u \in \mathbb{S}^n$ such that $K \subseteq \operatorname{int} \mathbb{S}_u^+$. Then*

$$H_{n-1}^{\mathbb{S}^n}(K, w) = H_{n-1}^{\mathbb{R}^n_u}(g_u(K), g_u(w)) \left( \frac{\cos(h_u(K, N_w^K | \mathbb{S}_u))}{\cos(d(u, w))} \right)^{n+1} \quad (4.13)$$

*for $\mathcal{H}^{n-1}$-almost all $w \in \operatorname{bd} K$.*

*Proof.* Since $g_u$ is a geodesic diffeomorphism we conclude that $H_{n-1}^{\mathbb{S}^n}(K, w) = 0$ if and only if $H_{n-1}^{\mathbb{R}^n_u}(g_u(K), g_u(w)) = 0$. Thus, in this case, (4.13) holds. So assume $H_{n-1}^{\mathbb{S}^n}(K, w) > 0$. By (4.9), we have $N_w^K | \mathbb{S}_u = N_{g_u(w)}^{g_u(K)}$, which implies that

$$R_{K^\circ} \circ N_w^K = N^{g_u(K)} \circ g_u(w),$$

where $N^{g_u(K)} \colon \operatorname{bd} g_u(K) \to \mathbb{S}_u$ is the Euclidean Gauss map of $g_u(K)$. For the outer unit normal $N^{g_u(K)}$ of a Euclidean convex body $L$ and a regular boundary point $x$ of $L$, we have $J^{\operatorname{bd} L}(N^L)(x) = H_{n-1}^{\mathbb{R}^n_u}(L, x)$ (see [**31**, Lemma 2.3]). By this, we conclude

$$J^{\mathbb{S}_u}(R_{K^\circ})(N_w^K) H_{n-1}^{\mathbb{S}^n}(K, w) = J^{\operatorname{bd} g_u(K)}(N^{g_u(K)})(g_u(w)) J^{\operatorname{bd} K}(g_u)(w),$$

and, by (4.12) and (4.10),

$$\frac{\sin(h_{-u}(K^\circ, N_{N_w^K}^{K^\circ} | \mathbb{S}_u))}{\sin(d(-u, N_w^K))^n} H_{n-1}^{\mathbb{S}^n}(K, w) = H_{n-1}^{\mathbb{R}^n_u}(g_u(K), g_u(w)) \frac{\cos(h_u(K, N_w^K | \mathbb{S}_u))}{\cos(d(u, w))^n}.$$

Since $H_{n-1}^{\mathbb{S}^n}(K, w) > 0$, we have that $N_w^K$ is the unique outer unit normal to $w$ and this also implies that $w$ is the unique outer unit normal to $N_w^K \in K^\circ$, thus, $(N^{K^\circ} \circ N^K)(w) = w$. By the duality formula (4.5), we obtain (4.13). We note that $\rho_{-u}(K^\circ, N_w^K | \mathbb{S}_u) = d(-u, N_w^K)$. q.e.d.

The following lemma can be considered as a spherical version of the splitting of Lebesgue integrals by orthogonal subspaces.

**Lemma 4.5** ([**66**]). *Let $S$ be a $k$-sphere, $0 \leq k \leq n-1$, and let $f \colon \mathbb{S}^n \to \mathbb{R}$ be a non-negative measurable function. Then*

$$\int_{\mathbb{S}^n} f(w)\, dw = \int_S \int_{\operatorname{conv}(S^\circ, v)} \sin(d(S^\circ, u))^k f(u)\, du\, dv.$$

For a Euclidean convex body $L \in \mathcal{K}_0(\mathbb{R}^{n+1})$ with $0 \in \operatorname{int} L$ we can express the volume of $L$ by integrating the cone volume measure over the boundary of $L$, i.e.,

$$\operatorname{vol}_n(L) = \frac{1}{n} \int_{\operatorname{bd} L} x \cdot N_x^L\, dx = \int_{\operatorname{bd} L} \frac{x \cdot N_x^L}{\|x\|^n} \int_0^{\|x\|} r^{n-1}\, dr\, dx.$$

The following proposition is a spherical version of this for $K \in \mathcal{K}_0^p(\mathbb{S}^n)$, where we fix a reference point $u \in \operatorname{int} K$.

**Proposition 4.6.** *Let $K \in \mathcal{K}_0^p(\mathbb{S}^n)$ and $u \in \operatorname{int} K$ such that $K \subseteq \operatorname{int} \mathbb{S}_u^+$. Then*

$$\operatorname{vol}_n(K) = \int_{\operatorname{bd} K} \frac{-u \cdot N_w^K}{\sin(d(u,w))^n} \int_0^{d(u,w)} \sin(t)^{n-1}\, dt\, dw.$$

*In particular, for $K, L \in \mathcal{K}_0^p(\mathbb{S}^n)$ such that $u \in \operatorname{int} L$ and $L \subseteq K \subseteq \operatorname{int} \mathbb{S}_u^+$,*

$$\operatorname{vol}_n(K) - \operatorname{vol}_n(L) = \int_{\operatorname{bd} K} \frac{-u \cdot N_w^K}{\sin(d(u,w))^n} \int_{d(u,w_L)}^{d(u,w)} \sin(t)^{n-1}\, dt\, dw,$$

*where we set $\{w_L\} = \operatorname{bd} L \cap \operatorname{conv}(u, w)$ for $w \in \operatorname{bd} K$.*

*Proof.* Using Lemma 4.5 and the fact that for the hypersphere $\mathbb{S}_u$ we have $\mathbb{S}_u^\circ = \{-u, u\}$, we obtain

$$\operatorname{vol}_n(K) = \int_{\mathbb{S}_u} \int_{\operatorname{conv}(\mathbb{S}_u^\circ, v)} \sin(d(z, \mathbb{S}_u^\circ))^{n-1} I_K(z)\, d\mathcal{H}^1(z)\, d\mathcal{H}^{n-1}(v). \qquad (*)$$

Now for any $z \in \operatorname{conv}(\mathbb{S}_u^\circ, v)$ we can write $z = \cos(t)u + \sin(t)v$ where $t$ is determined by $t = d(\mathbb{S}_u^\circ, z) = d(u, z)$. From definition (4.4) of $\rho_u(K, v)$, we conclude

$$(*) = \int_{\mathbb{S}_u} \int_0^{\rho_u(K,v)} \sin(t)^{n-1}\, dt\, dv = \int_{\operatorname{bd} K} J^{\operatorname{bd} K}(R_K)(w) \int_0^{\rho_u(K, R_K(w))} \sin(t)^{n-1}\, dt\, dw,$$

where we used the area formula for the spherical radial map $R_K \colon \operatorname{bd} K \to \mathbb{S}_u$ defined by $R_K(w) = w|\mathbb{S}_u$. By (4.12) and since, for any $w \in \operatorname{bd} K$, we have $\rho_u(K, R_K(w)) = d(u,w)$, we are done.

The second statement of the proposition follows easily. q.e.d.

## 5. The Spherical Convex Floating Body

In this section we collect results about the spherical convex floating body which we will need in Section 6 to prove Theorem 2.1.

**Definition** (Spherical Convex Floating Body). *Let $K \in \mathcal{K}_0(\mathbb{S}^n)$ and $\delta > 0$. Unless it is empty, the (spherical) convex floating body $K_{[\delta]}$ is defined by*

$$K_{[\delta]} = \bigcap \{\mathbb{S}_z^- : z \in \mathbb{S}^n \text{ such that } \mathrm{vol}_n(K \cap \mathbb{S}_z^+) \leq \delta\}. \tag{5.1}$$

$K_{[\delta]}$ is convex because it is an intersection of closed hemispheres and, as we will show, it exists if $\delta$ is small enough. First we show that for a proper convex body the intersection can be parametrized with respect to a fixed hypersphere.

**Lemma 5.1.** *Let $K \in \mathcal{K}_0^p(\mathbb{S}^n)$, $\delta \in (0, \mathrm{vol}_n(K))$ and $u \in \mathrm{int}\, K$ such that $K \subseteq \mathrm{int}\, \mathbb{S}_u^+$. For $v \in \mathbb{S}_u$ there exists a unique $s(v, \delta) \in (-\frac{\pi}{2}, \frac{\pi}{2})$ determined by*

$$\mathrm{vol}_n\left(K \cap \mathbb{S}_{u,v,s(v,\delta)}^+\right) = \delta,$$

*where $\mathbb{S}_{u,v,s(v,\delta)}^+ = \{w \in \mathbb{S}^n : v \cdot w \geq \tan(s(v,\delta))(u \cdot w)\}$. Moreover,*

$$K_{[\delta]} = \bigcap_{v \in \mathbb{S}_u} \mathbb{S}_{u,v,s(v,\delta)}^-, \tag{5.2}$$

*and $s(v, \delta)$ is continuous in both arguments and strictly decreasing in $\delta$.*

*Proof.* We consider the function $f$ defined by $f(v, s) = \mathrm{vol}_n(K \cap \mathbb{S}_{u,v,s}^+)$. Then $f$ is continuous in both arguments and strictly decreasing in $s$. By (4.2), we have, for $z = \cos(h_u(K, v))v - \sin(h_u(K, v))u$, that $\mathbb{S}_z$ is a supporting hypersphere to $K$. Therefore $f(v, h_u(K, v)) = 0$. Similarly, we have that $f(v, -h_u(K, -v)) = \mathrm{vol}_n(K)$. We therefore conclude that there exists a unique $s(v, \delta) \in (-h_u(K, -v), h_u(K, v)) \subseteq (-\frac{\pi}{2}, \frac{\pi}{2})$ such that $f(v, s(v, \delta)) = \delta$. Thus $s(v, \delta)$ is continuous in both arguments and strictly decreasing in $\delta$.

To prove (5.2) we only have to show that

$$K_{[\delta]} \supseteq \bigcap_{v \in \mathbb{S}_u} \mathbb{S}_{u,v,s(v,\delta)}^-. \tag{5.3}$$

Let $z \in \mathbb{S}^n$ such that $\mathrm{vol}_n(K \cap \mathbb{S}_z^+) \leq \delta$. For $z \neq \pm u$, we set $v = (z|\mathbb{S}_u)$. Then there is a unique $s' \in (-\frac{\pi}{2}, \frac{\pi}{2})$ such that $z = \cos(s')v - \sin(s')u$. We conclude

$$f(s') = \mathrm{vol}_n\left(K \cap \mathbb{S}_{u,v,s'}^+\right) = \mathrm{vol}_n(K \cap \mathbb{S}_z^+) \leq \delta = f(s(v, \delta)).$$

Thus, $s' \geq s(v, \delta)$, and therefore,

$$\mathbb{S}_u^+ \cap \mathbb{S}_z^- = \mathbb{S}_u^+ \cap \mathbb{S}_{u,v,s'}^- = \{w \in \mathbb{S}^n : u \cdot w \geq 0 \text{ and } v \cdot w \leq \tan(s')(u \cdot w)\}$$

$$\supseteq \mathbb{S}_u^+ \cap \mathbb{S}_{u,v,s(v,\delta)}^-.$$

Since $\mathbb{S}_u^+ \supseteq K \supseteq K_{[\delta]}$, we obtain

$$K_{[\delta]} = \mathbb{S}_u^+ \cap K_{[\delta]} = \bigcap\{\mathbb{S}_u^+ \cap \mathbb{S}_z^- : \mathrm{vol}_n(K \cap \mathbb{S}_z^+) \leq \delta\}$$

$$\supseteq \bigcap\{\mathbb{S}_u^+ \cap \mathbb{S}_{u,v,s(v,\delta)}^- : v \in \mathbb{S}_u\} = \mathbb{S}_u^+ \cap \left(\bigcap_{v \in \mathbb{S}_u} \mathbb{S}_{u,v,s(v,\delta)}^-\right).$$

Since $\mathbb{S}_{u,v,s(v,\delta)}^- \cap \mathbb{S}_{u,-v,s(-v,\delta)}^- \subseteq \mathbb{S}_u^+$, we conclude (5.3). q.e.d.

The following theorem relates the Euclidean floating body of the gnomonic projection of a proper spherical convex body to the spherical convex floating body.

**Theorem 5.2.** *Let $K \in \mathcal{K}_0^p(\mathbb{S}^n)$ such that $C_u(\alpha) \subseteq \mathrm{int}\, K$ and $K \subseteq C_u(\beta)$ for some $u \in \mathbb{S}^n$ and $\alpha, \beta \in [0, \frac{\pi}{2})$. Then, for $\delta > 0$ small enough, we have*

$$g_u(K)_{\left[\frac{\delta}{\cos(\beta)^{n+1}}\right]} \subseteq g_u(K_{[\delta]}) \subseteq g_u(K)_{\left[\frac{\delta}{\cos(\alpha)^{n+1}}\right]}.$$

In particular, this shows that $K_{[\delta]}$ exists if $\delta$ is small enough.

*Proof.* Set $L = g_u(K)$. By (5.2), we have $K_{[\delta]} = \bigcap_{v \in \mathbb{S}_u} \mathbb{S}^-_{u,v,s(v,\delta)}$, where $s(v,\delta)$ is determined by $\delta = \text{vol}_n\left(K \cap \mathbb{S}^+_{u,v,s(v,\delta)}\right)$. The gnomonic projection $g_u$ maps $\mathbb{S}_{u,v,s(v,\delta)}$ to the hyperplane $\mathbb{H}_{v,\tan(s(v,\delta))} = \{x \in \mathbb{R}^n_u : x \cdot v = \tan(s(v,\delta))\}$. For $\delta$ small enough, we have for $v \in \mathbb{S}_u$ that

$$\emptyset = g_u\left(\left(K \cap \mathbb{S}_{u,v,s(v,\delta)}\right) \cap C_u(\alpha)\right) = \left(L \cap \mathbb{H}^+_{v,\tan(s(v,\delta))}\right) \cap B(0, \tan(\alpha)).$$

By (4.8), we conclude

$$\delta = \text{vol}_n\left(K \cap \mathbb{S}^+_{u,v,s(v,\delta)}\right) = \int_{L \cap \mathbb{H}^+_{v,\tan(s(v,\delta))}} \frac{dx}{(1+\|x\|^2)^{\frac{n+1}{2}}}$$

$$\leq \cos(\alpha)^{n+1} \text{vol}_n\left(L \cap \mathbb{H}^+_{v,\tan(s(v,\delta))}\right).$$

Now let $\widetilde{s}(v,\delta)$ such that $\delta = \cos(\alpha)^{n+1} \text{vol}_n\left(L \cap \mathbb{H}^+_{v,\widetilde{s}(v,\delta)}\right)$. Then obviously $\widetilde{s}(v,\delta) \geq \tan(s(v,\delta))$ and $\mathbb{H}^-_{v,\tan(s(v,\delta))} \subseteq \mathbb{H}^-_{v,\widetilde{s}(v,\delta)}$. By (3.1) we deduce

$$g_u(K_{[\delta]}) = \bigcap_{v \in \mathbb{S}_u} \mathbb{H}^-_{v,\tan(s(v,\delta))} \subseteq \bigcap_{v \in \mathbb{S}_u} \mathbb{H}^-_{v,\widetilde{s}(v,\delta)} = L_{\left[\frac{\delta}{\cos(\alpha)^{n+1}}\right]}.$$

For the converse, since $g_u(K) \subseteq g_u(C_u(\beta)) = B(0, \tan(\beta))$, we first note that

$$\delta = \text{vol}_n\left(K \cap \mathbb{S}^+_{u,v,s(v,\delta)}\right) = \int_{L \cap \mathbb{H}^+_{v,\tan(s(v,\delta))}} \frac{dx}{(1+\|x\|^2)^{\frac{n+1}{2}}}$$

$$\geq \cos(\beta)^{n+1} \text{vol}_n\left(L \cap \mathbb{H}^+_{\tan(s(v,\delta))}\right).$$

Now let $\widetilde{S}(v,\delta)$ such that $\delta = \cos(\beta)^{n+1} \text{vol}_n\left(L \cap \mathbb{H}^+_{v,\widetilde{S}(v,\delta)}\right)$. Then we have $\widetilde{S}(v,\delta) \leq \tan(s(v,\delta))$ and therefore $\mathbb{H}^-_{v,\tan(s(v,\delta))} \supseteq \mathbb{H}^-_{v,\widetilde{S}(v,\delta)}$. We conclude

$$g_u(K_{[\delta]}) = \bigcap_{v \in \mathbb{S}_u} \mathbb{H}^-_{v,\tan(s(v,\delta))} \supseteq \bigcap_{v \in \mathbb{S}_u} \mathbb{H}^-_{v,\widetilde{S}(v,\delta)} = L_{\left[\frac{\delta}{\cos(\beta)^{n+1}}\right]}.$$

q.e.d.

In the following three lemmas we establish properties of the spherical convex floating body as $\delta$ goes to 0. First we show that the boundary points of the floating body converge to boundary points of the convex body.

**Lemma 5.3.** *If $K \in \mathcal{K}_0(\mathbb{S}^n)$ and $\delta_1 \leq \delta_2$, then $K_{[\delta_1]} \supseteq K_{[\delta_2]}$. In particular, we have* $\text{int } K = \bigcup_{\delta > 0} K_{[\delta]}$.

*Furthermore, let $K \in \mathcal{K}^p_0(\mathbb{S}^n)$ and $u \in \text{int } K$ such that $K \subseteq \text{int } \mathbb{S}^+_u$. For $w \in \text{bd } K$, put $\{w_\delta\} = \text{bd } K_{[\delta]} \cap \text{conv}(u, w)$. Then $\lim_{\delta \to 0} w_\delta = w$.*

*Proof.* The monotonicity of the floating body is obvious from its definition.

First we prove $\text{int } K = \bigcup_{\delta > 0} K_{[\delta]}$. Let $\delta > 0$ small enough such that $K_{[\delta]}$ is non-empty. If $w \in K_{[\delta]}$, then every hypersphere through $w$ cuts off a set of $K$ of volume greater or equal to $\delta$. Thus $K_{[\delta]} \subseteq \text{int } K$ and we conclude $\bigcup_{\delta > 0} K_{[\delta]} \subseteq \text{int } K$.

In order to prove the converse, we assume that there is $w \in \operatorname{int} K$ such that $w \notin \bigcap_{\delta > 0} K_{[\delta]}$. For every $v \in \mathbb{S}_w$ we have

$$\operatorname{vol}_n\left(K \cap \mathbb{S}_v^+\right) > 0. \tag{5.4}$$

Since $w \notin \bigcap_{\delta > 0} K_{[\delta]}$ we conclude that for every $\delta > 0$ there exists $v(\delta) \in \mathbb{S}_w$ such that $\operatorname{vol}_n\left(K \cap \mathbb{S}_{v(\delta)}^+\right) < \delta$. By compactness of $\mathbb{S}_w$ and continuity there exists $v_0 \in \mathbb{S}_w$ such that $\operatorname{vol}_n\left(K \cap \mathbb{S}_{v_0}^+\right) = 0$. This is a contradiction to (5.4).

Finally, let $K \in \mathcal{K}_0^p(\mathbb{S}^n)$ and $u \in \operatorname{int} K$ such that $K \subseteq \operatorname{int} \mathbb{S}_u^+$. We have

$$\bigcup_{\delta > 0} \operatorname{conv}(u, w_\delta) = \bigcup_{\delta > 0} K_{[\delta]} \cap \operatorname{conv}(u, w) = \operatorname{conv}(u, w) \backslash \{w\}.$$

We conclude $\lim_{\delta \to 0} d(w_\delta, w) = 0$. \hfill q.e.d.

In the next lemma we show that the outer normals of the spherical convex floating body converge to the outer normals of the convex body.

**Lemma 5.4.** *Let $K \in \mathcal{K}_0^p(\mathbb{S}^n)$, $u \in \operatorname{int} K$ and $w \in \operatorname{bd} K$ be a regular boundary point. For $\delta > 0$ such that $K_{[\delta]} \neq \emptyset$, we set $\{w_\delta\} = \operatorname{bd} K \cap \operatorname{conv}(u, w)$. Then*

$$\lim_{\delta \to 0^+} N_{w_\delta}^{K_{[\delta]}} = N_w^K, \tag{5.5}$$

*where $N_{w_\delta}^{K_{[\delta]}}$ is an outer unit normal to $K_{[\delta]}$ in $w_\delta$ such that*

$$\delta = \operatorname{vol}_n\left(K \cap \mathbb{S}_{N_{w_\delta}^{K_{[\delta]}}}^+\right).$$

*In particular, for all $\varepsilon > 0$ there exists $\delta(\varepsilon)$ such that*

$$N_w^K \cdot N_{w_\delta}^{K_{[\delta]}} \geq 1 - \varepsilon, \tag{5.6}$$

*for all $\delta \leq \delta(\varepsilon)$ and, if $K \subseteq \operatorname{int} \mathbb{S}_u^+$, then for all $\delta \leq \delta(\varepsilon)$*

$$(N_w^K | \mathbb{S}_u) \cdot (N_{w_\delta}^{K_{[\delta]}} | \mathbb{S}_u) \geq 1 - \varepsilon. \tag{5.7}$$

*Proof.* Suppose (5.5) is not true. Then, by compactness, there exists a subsequence $(\delta_i)_{i \in \mathbb{N}}$ such that $\lim_{i \to \infty} \delta_i = 0$, $\lim_{i \to \infty} N_{w_{\delta_i}}^{K_{[\delta_i]}} = v_0$ and $v_0 \neq N_w^K$. By the choice of $N_{w_{\delta_i}}^{K_{[\delta_i]}}$, we have $\operatorname{vol}_n\left(K \cap \mathbb{S}_{N_{w_{\delta_i}}^{K_{[\delta_i]}}}^+\right) = \delta_i$. We conclude that $\operatorname{vol}_n(K \cap \mathbb{S}_{v_0}^+) = 0$. By Lemma 5.3, we have $\lim_{i \to \infty} w_{\delta_i} = w$, thus $v_0$ is a normal to $\operatorname{bd} K$ in $w$. This contradicts the assumption that $w$ is a regular boundary point and therefore has a unique outer unit normal $N_w^K \neq v_0$.

The other statements, (5.6) and (5.7), follow directly from (5.5). \hfill q.e.d.

Let $w \in \operatorname{bd} K$ be a boundary point such that $H_{n-1}^{\mathbb{S}^n}(K, w) > 0$. Then $K \cap \mathbb{S}_{w, N_w^K, -\Delta}^+$ is decreasing in $\Delta$ and for $\Delta = 0$ equals just $\{w\}$. Thus for $\Delta$ small enough, $K \cap \mathbb{S}_{w, N_w^K, -\Delta}^+$ is contained in some arbitrarily small cap around $w$. By continuity this will still be true if we allow directions $v \in \mathbb{S}_w$ close to $N_w^K$. Let $w_\delta$ be a boundary point of the floating body $K_{[\delta]}$ that converges to a boundary point $w \in \operatorname{bd} K$. Then, by Lemma 5.4, the normals to $w_\delta$ converge to $N_w^K$. Thus, if we consider directions $v \in \mathbb{S}_w$ close to $N_w^K$, then, for $\delta$ small enough, $w_\delta$ will already be determined by these directions. Hence, if we replace $K$ by $K' = K \cap C_w(\varepsilon)$ for arbitrarily small $\varepsilon$, then for $\delta$ small enough, we will have $w_\delta^{K'} = w_\delta^K$. We will prove this rigorously in the following lemma.

**Lemma 5.5.** *Let $K \in \mathcal{K}_0^p(\mathbb{S}^n)$ and $w \in \operatorname{bd} K$ such that $H_{n-1}^{\mathbb{S}^n}(K,w) > 0$. For $\varepsilon > 0$ set $K' = K \cap C_w(\varepsilon)$.*

*(i) There exists $\Delta_\varepsilon$ such that for all $\Delta < \Delta_\varepsilon$ we have*
$$K' \cap \mathbb{S}^+_{w,N_w^K,-\Delta} = K \cap \mathbb{S}^+_{w,N_w^K,-\Delta}.$$

*(ii) There exists $\xi_\varepsilon$ and $\eta_\varepsilon$ such that, for all $v \in \mathbb{S}_w$ with $d(v, N_w^K) < \xi_\varepsilon$ and $\Delta < \eta_\varepsilon$, we have*
$$K' \cap \mathbb{S}^+_{w,v,-\Delta} = K \cap \mathbb{S}^+_{w,v,-\Delta}.$$

*(iii) Let $u \in \operatorname{int} K'$. There exists $\delta_\varepsilon$ such that, for all $\delta < \delta_\varepsilon$, we have*
$$K'_{[\delta]} \cap \operatorname{conv}(u,w) = K_{[\delta]} \cap \operatorname{conv}(u,w).$$

*In particular, we have $w_\delta^{K'} = w_\delta^K$ for all $\delta < \delta_\varepsilon$.*

*Proof.* (i) Assume this does not hold. Then there exists an $\varepsilon > 0$ such that for all $\Delta > 0$, we have
$$\emptyset \neq K \cap \mathbb{S}^+_{w,N_w^K,-\Delta} \setminus \left(K' \cap \mathbb{S}^+_{w,N_w^K,-\Delta}\right)$$
$$= (K \setminus C_w(\varepsilon)) \cap \mathbb{S}^+_{w,N_w^K,-\Delta} \subseteq (K \setminus \operatorname{int} C_w(\varepsilon)) \cap \mathbb{S}^+_{w,N_w^K,-\Delta}.$$

For $\Delta_1 \leq \Delta_2$, we have $(K \setminus \operatorname{int} C_w(\varepsilon)) \cap \mathbb{S}^+_{w,N_w^K,-\Delta_1} \subseteq (K \setminus \operatorname{int} C_w(\varepsilon)) \cap \mathbb{S}^+_{w,N_w^K,-\Delta_2}$, and, by compactness, we conclude that $\emptyset \neq (K \setminus \operatorname{int} C_w(\varepsilon)) \cap \mathbb{S}_{N_w^K}$. Since $\mathbb{S}_{N_w^K}$ is a supporting hyperplane at $w$, this implies that there exists $v \in K \cap \mathbb{S}_{N_w^K}$ such that $d(v, w) \geq \varepsilon$. Since $K$ is convex, the whole segment $\operatorname{conv}(v, w)$ is contained in $K$. Considering $L = g_w(K \cap C_w(\varepsilon))$, we see that the boundary of $L$ contains the segment $g_w(C_w(\varepsilon) \cap \operatorname{conv}(w,v))$, and we conclude that $H_{n-1}^{\mathbb{R}_w^n}(L,0) = 0$. This implies, by Lemma 4.4, that $H_{n-1}^{\mathbb{S}^n}(K,w) = 0$ which is a contradiction.

*(ii)* We may assume $\varepsilon < \frac{\pi}{3}$. By (i), there exists $\Delta_{\varepsilon/2}$ such that $\Delta_{\varepsilon/2} < \frac{\varepsilon}{2}$ and
$$(K \setminus C_w(\varepsilon/2)) \cap \mathbb{S}^+_{w,N_w^K,-\Delta_{\varepsilon/2}} = \emptyset. \tag{5.8}$$

We set $\eta_\varepsilon = \frac{\Delta_{\varepsilon/2}}{2} < \frac{\varepsilon}{4}$ and
$$\sin(\xi_\varepsilon) = \frac{\tan(\Delta_{\varepsilon/2}) - \tan(\eta_\varepsilon)}{\sqrt{\tan(\varepsilon)^2 - \tan(\eta_\varepsilon)^2}} > 0.$$

We show that, for all $v \in \mathbb{S}_w$ such that $d(v, N_w^K) < \xi_\varepsilon$, we have
$$(K \setminus C_w(\varepsilon)) \cap \mathbb{S}^+_{w,v,-\eta_\varepsilon} = \emptyset. \tag{5.9}$$

This implies $K' \cap \mathbb{S}^+_{u,v,-\Delta} = K \cap \mathbb{S}^+_{u,v,-\Delta}$ for all $\Delta \leq \eta_\varepsilon$.

Assume (5.9) is not true. Then there exists $z \in (K \setminus C_w(\varepsilon)) \cap \mathbb{S}^+_{w,v,-\eta_\varepsilon}$. Since $K$ is convex, the whole segment $\operatorname{conv}(z,w)$ is contained in $K$. Since $d(z,w) > \varepsilon$, there exists $\{z'\} = \operatorname{bd} C_w(\varepsilon) \cap \operatorname{conv}(z,w)$. We will show that $z' \in \mathbb{S}^+_{w,N_w^K,-\Delta_{\varepsilon/2}}$. This will be a contradiction to (5.8), since by construction $z' \in K \setminus C_w(\varepsilon/2)$. We have to show that
$$\frac{-z' \cdot N_w^K}{z' \cdot w} \leq \tan(\Delta_{\varepsilon/2}). \tag{5.10}$$

Since $z' \in \operatorname{bd} C_w(\varepsilon)$, we have $z' = \cos(\varepsilon)w + \sin(\varepsilon)(z'|\mathbb{S}_w)$ and conclude
$$\frac{-z' \cdot N_w^K}{z' \cdot w} = \tan(\varepsilon) \sin\left(d(z'|\mathbb{S}_w, N_w^K) - \frac{\pi}{2}\right),$$

and, since $z' \in \mathbb{S}^+_{w,v,-\eta_\varepsilon}$, we obtain

$$\tan(\eta_\varepsilon) \geq \frac{-z' \cdot v}{z' \cdot w} = \tan(\varepsilon) \sin\left(d(z'|\mathbb{S}_w, v) - \frac{\pi}{2}\right).$$

Thus $d(z'|\mathbb{S}_w, v) \leq \frac{\pi}{2} + \arcsin\left(\frac{\tan(\eta_\varepsilon)}{\tan(\varepsilon)}\right)$ and by the triangle inequality we have

$$d(z'|\mathbb{S}_w, N_w^K) \leq d(z'|\mathbb{S}_w, v) + d(v, N_w^K) \leq \frac{\pi}{2} + \arcsin\left(\frac{\tan(\eta_\varepsilon)}{\tan(\varepsilon)}\right) + \xi_\varepsilon.$$

Finally, we obtain (5.10) from

$$\frac{-z' \cdot N_w^K}{z' \cdot w} \leq \tan(\varepsilon) \sin\left(\arcsin\left(\frac{\tan(\eta_\varepsilon)}{\tan(\varepsilon)}\right) + \xi_\varepsilon\right)$$
$$= \sin(\xi_\varepsilon)\sqrt{\tan(\varepsilon)^2 - \tan(\eta_\varepsilon)^2} + \cos(\xi_\varepsilon)\tan(\eta_\varepsilon)$$
$$\leq \sin(\xi_\varepsilon)\sqrt{\tan(\varepsilon)^2 - \tan(\eta_\varepsilon)^2} + \tan(\eta_\varepsilon) \leq \tan(\Delta_{\varepsilon/2}).$$

(iii) Since $K' \subseteq C_w(\varepsilon)$, by (5.2), we can write $K'_{[\delta]} = \bigcap_{v \in \mathbb{S}_w} \mathbb{S}^-_{w,v,s^{K'}(v,\delta)}$. Here $s^{K'}(v,\delta)$ is uniquely determined by $\delta = \text{vol}_n\left(K' \cap \mathbb{S}^+_{w,v,s^{K'}(v,\delta)}\right)$ and is continuous in both arguments. By (ii), there exists $\xi_\varepsilon$ and $\eta_\varepsilon$ such that, for all $v \in \mathbb{S}_w$ with $d(v, N_w^K) < \xi_\varepsilon$, we have

$$\delta^{K'}(v, -\Delta) = (\text{vol}_n(K' \cap \mathbb{S}_{w,v,-\Delta}) = \text{vol}_n(K \cap \mathbb{S}_{w,v,-\Delta}) = \delta^K(v, -\Delta),$$

for all $\Delta < \eta_\varepsilon$. Hence $s^K(v, \delta)$ exists for $v \in \mathbb{S}_w$ such that $d(v, N_w^K) < \xi_\varepsilon$ and $\delta$ small. Thus, there exist $\delta_1$ such that for all $\delta < \delta_1$ and $v \in \mathbb{S}_w$ such that $d(v, N_w^{K'}) = d(v, N_w^K) < \xi_\varepsilon$, we have $s^{K'}(v, \delta) = s^K(v, \delta)$.

**Claim:** There exists $\delta_2$ such that, for all $\delta < \delta_2$, we have

$$K_{[\delta]} \cap \text{conv}(u, w) = \bigcap \{\mathbb{S}^-_{w,v,s^K(v,\delta)} : v \in \mathbb{S}_w, d(v, N_w^K) \leq \xi_\varepsilon\} \cap \text{conv}(u, w).$$

Assume that this is not true. Then, for all $\delta > 0$, we have

$$\text{conv}(u, w_\delta^K) = K_{[\delta]} \cap \text{conv}(u, w)$$
$$\subsetneq \bigcap \{\mathbb{S}^-_{w,v,s^K(v,\delta)} : d(v, N_w^K) \leq \xi_\varepsilon\} \cap \text{conv}(u, w).$$

Thus, for any normal $z$ to $K_{[\delta]}$ in $w_\delta^K$, we have $d(z|\mathbb{S}_w, N_w^K) \geq \xi_\varepsilon > 0$ for all $\delta > 0$. This is a contradiction to Lemma 5.4.

With a similar argument for $K'$ we also find a $\delta_3$ such that, for all $\delta < \delta_3$,

$$K'_{[\delta]} \cap \text{conv}(u, w) = \bigcap \{\mathbb{S}^-_{w,v,s^{K'}(v,\delta)} : v \in \mathbb{S}_w, d(v, N_w^{K'}) \leq \xi_\varepsilon\} \cap \text{conv}(u, w).$$

Setting $\delta_\varepsilon = \min\{\delta_1, \delta_2, \delta_3\}$, we conclude, for all $\delta < \delta_\varepsilon$,

$$K'_{[\delta]} \cap \text{conv}(u, w) = \bigcap \{\mathbb{S}^-_{w,v,-s^{K'}(v,\delta)} : d(v, N_w^{K'}) < \xi_1\} \cap \text{conv}(u, w)$$
$$= \bigcap \{\mathbb{S}^-_{w,v,-s^K(v,\delta)} : d(v, N_w^K) < \xi_1\} \cap \text{conv}(u, w)$$
$$= K_{[\delta]} \cap \text{conv}(u, w).$$

The second statement of (iii) is obvious since $K_{[\delta]} \cap \text{conv}(u, w) = \text{conv}(u, w_\delta^K)$ and $K'_{[\delta]} \cap \text{conv}(u, w) = \text{conv}(u, w_\delta^{K'})$. q.e.d.

## 6. Proof of the Main Result

We are now ready to proof Theorem 2.1. By Lemma 4.1, there exists $u \in \operatorname{int} K$ such that $K \subseteq \operatorname{int} \mathbb{S}_u^+$. Using Proposition 4.6, we may write the left hand side of (2.3) as

$$\frac{\operatorname{vol}_n(K) - \operatorname{vol}_n(K_{[\delta]})}{\delta^{\frac{2}{n+1}}} = \int_{\operatorname{bd} K} \underbrace{\delta^{-\frac{2}{n+1}} \frac{-u \cdot N_w^K}{\sin(d(w,u))^n} \int_{d(w_\delta, u)}^{d(w,u)} \sin(t)^{n-1} \, dt}_{f(w,\delta)} \, dw. \tag{6.1}$$

The proof of Theorem 2.1 will now be completed in two steps. We will first show that the integrand $f(w,\delta)$ in the integral on the right hand side of (6.1) is bounded from above uniformly in $\delta$ almost everywhere by an integrable function.

**Lemma 6.1.** *Let $K \in \mathcal{K}_0^p(\mathbb{S}^n)$ and $u \in \operatorname{int} K$ such that $K \subseteq \operatorname{int} \mathbb{S}_u^+$. There exists $\delta_0 > 0$ and an integrable function $g \colon \operatorname{bd} K \to \mathbb{R}$ such that, for all $\delta < \delta_0$,*

$$\delta^{-\frac{2}{n+1}} \frac{-u \cdot N_w^K}{\sin(d(w,u))^n} \int_{d(w_\delta, u)}^{d(w,u)} \sin(t)^{n-1} \, dt \leq g(w) \tag{6.2}$$

*for $\mathcal{H}^{n-1}$-almost every $w \in \operatorname{bd} K$, where $\{w_\delta\} = \operatorname{bd} K_{[\delta]} \cap \operatorname{conv}(u, w)$.*

*Proof.* Since $u \in \operatorname{int} K$ and $K \subseteq \operatorname{int} \mathbb{S}_u^+$, there is $0 < \alpha < \beta < \frac{\pi}{2}$ such that $C_u(\alpha) \subseteq \operatorname{int} K$ and $K \subseteq C_u(\beta)$. Therefore, $\sin(d(u,w)) \geq \sin(\alpha)$, $-u \cdot N_w^K \leq 1$ and conclude

$$\delta^{-\frac{2}{n+1}} \frac{-u \cdot N_w^K}{\sin(d(w,u))^n} \int_{d(w_\delta, u)}^{d(w,u)} \sin(t)^{n-1} \, dt \leq \frac{1}{\sin(\alpha)} \frac{d(w, w_\delta)}{\delta^{\frac{2}{n+1}}}.$$

We will show that there exists $C > 0$ and $\delta_0$ such that, for all $\delta < \delta_0$,

$$\frac{1}{\sin(\alpha)} \frac{d(w, w_\delta)}{\delta^{\frac{2}{n+1}}} \leq C r_{g_u(K)}(g_u(w))^{-\frac{n-1}{n+1}} \tag{6.3}$$

for $\mathcal{H}^{n-1}$-almost all $w \in \operatorname{bd} K$. Then the right hand side of (6.3) is integrable: By (4.11), the fact that $1 + (x \cdot N_x^{g_u(K)})^2 \leq 1 + \|x\|^2$ and $1 + \|x\|^2 \geq 1$, we have

$$\int_{\operatorname{bd} K} r_{g_u(K)}(g_u(w))^{-\frac{n-1}{n+1}} \, dw = \int_{\operatorname{bd} g_u(K)} r_{g_u(K)}(x)^{-\frac{n-1}{n+1}} \frac{\sqrt{1 + (x \cdot N_x^{g_u(K)})^2}}{(1 + \|x\|^2)^{\frac{n}{2}}} \, dx$$

$$\leq \int_{\operatorname{bd} g_u(K)} r_{g_u(K)}(x)^{-\frac{n-1}{n+1}} \, dx,$$

which is finite by Theorem 3.3.

In order to prove (6.3), we set $L = g_u(K)$, $x = g_u(w)$. By (4.3), we have $\tan(d(u,w)) = \|g_u(w)\|$ and $\tan(d(u,w_\delta)) = \|g_u(w_\delta)\|$. We derive

$$\frac{d(w, w_\delta)}{\delta^{\frac{2}{n+1}}} = \frac{\arctan(\|x\|) - \arctan(\|g_u(w_\delta)\|)}{\delta^{\frac{2}{n+1}}} \leq \frac{\|x - g_u(w_\delta)\|}{\delta^{\frac{2}{n+1}}}.$$

By Theorem 5.2 and as $g_u(w_\delta)$ is on the line $g_u(\mathrm{conv}(u,w)) = \mathrm{conv}(0,x)$, we have $g_u(K_{[\delta]}) \subseteq L_{[\tilde\delta]}$, where $\tilde\delta = \frac{\delta}{\cos(\alpha)^{n+1}}$. Setting $\{x_{\tilde\delta}\} = \mathrm{bd}\, L_{[\tilde\delta]} \cap \mathrm{conv}(0,x)$, we conclude that $\|g_u(w_\delta)\| \geq \|x_{\tilde\delta}\|$. Therefore,

$$\frac{1}{\sin(\alpha)} \frac{d(w,w_\delta)}{\delta^{\frac{2}{n+1}}} \leq \frac{1}{\cos(\alpha)^2 \sin(\alpha)} \frac{\|x - x_{\tilde\delta}\|}{\tilde\delta^{\frac{2}{n+1}}}.$$

By Theorem 3.5, there exists $\tilde\delta_0$ and $\tilde C > 0$ such that, for all $\tilde\delta < \tilde\delta_0$,

$$\frac{\|x - x_{\tilde\delta}\|}{\tilde\delta^{\frac{2}{n+1}}} \leq \tilde C r_L(x)^{-\frac{n-1}{n+1}},$$

for $\mathcal{H}^{n-1}$-almost all $x \in \mathrm{bd}\, L$. Thus, for all $\delta < \tilde\delta_0 \cos(\alpha)^{n+1}$, we have

$$\frac{1}{\sin(\alpha)} \frac{d(w,w_\delta)}{\delta^{\frac{2}{n+1}}} \leq \frac{\tilde C}{\cos(\alpha)^2 \sin(\alpha)} r_{g_u(K)}(g_u(w))^{-\frac{n-1}{n+1}}$$

for $\mathcal{H}^{n-1}$-almost all $w \in \mathrm{bd}\, K$. q.e.d.

Using this lemma and Dominated Convergence Theorem, we can express the limit as $\delta$ tends to zero of the right hand side of (6.1) by the integral over the point-wise limit of the integrand.

**Theorem 6.2.** *Let $K \in \mathcal{K}_0^p(\mathbb{S}^n)$ and $u \in \mathrm{int}\, K$ such that $K \subseteq \mathrm{int}\, \mathbb{S}_u^+$. Then, for $\mathcal{H}^{n-1}$-almost all $w \in \mathrm{bd}\, K$, we have*

$$\lim_{\delta \to 0^+} \delta^{-\frac{2}{n+1}} \frac{-u \cdot N_w^K}{\sin(d(u,w))^n} \int_{d(u,w_\delta)}^{d(u,w)} \sin(t)^{n-1}\, dt = c_n H_{n-1}^{\mathbb{S}^n}(K,w)^{\frac{1}{n+1}}, \quad (6.4)$$

*where $\{w_\delta\} = \mathrm{bd}\, K_{[\delta]} \cap \mathrm{conv}(u,w)$ and $c_n = \frac{1}{2}\left(\frac{n+1}{\kappa_{n-1}}\right)^{\frac{2}{n+1}}$.*

*Proof.* For $t \in [d(u,w_\delta), d(u,w)]$ we have

$$\frac{\sin(d(u,w_\delta))}{\sin(d(u,w))} \leq \frac{\sin(t)}{\sin(d(u,w))} \leq 1.$$

Furthermore, $\lim_{\delta \to 0} d(u,w_\delta) = d(u,w)$ and $d(u,w) - d(u,w_\delta) = d(w,w_\delta)$. We therefore conclude

$$\lim_{\delta \to 0^+} \frac{-u \cdot N_w^K}{\sin(d(u,w))^n} \frac{1}{\delta^{\frac{2}{n+1}}} \int_{d(u,w_\delta)}^{d(u,w)} \sin(t)^{n-1} = \lim_{\delta \to 0^+} \frac{-u \cdot N_w^K}{\sin(d(u,w))} \frac{d(w,w_\delta)}{\delta^{\frac{2}{n+1}}}. \quad (6.5)$$

We first show that, for a regular boundary point $w$ of $K$ with positive curvature, the right hand side of (6.5) only depends on the local structure of $\mathrm{bd}\, K$ at $w$.

***Claim:*** *Let $\varepsilon \in \left(0, \frac{\pi}{2}\right)$ and $w \in \mathrm{bd}\, K$ such that $H_{n-1}^{\mathbb{S}^n}(K,w) > 0$. If we set $K' = K \cap C_w(\varepsilon)$ and let $u' \in \mathrm{int}\, K' \cap \mathrm{conv}(u,w)$ be such that $K' \subseteq \mathrm{int}\, \mathbb{S}_{u'}^+$, then, for $\delta$ small enough, we have*

$$\frac{-u \cdot N_w^K}{\sin(d(u,w))} \frac{d(w, w_\delta^K)}{\delta^{\frac{2}{n+1}}} = \frac{-u' \cdot N_w^K}{\sin(d(u',w))} \frac{d(w, w_\delta^{K'})}{\delta^{\frac{2}{n+1}}}. \quad (6.6)$$

Since $u'|\mathbb{S}_w = \frac{u' - \cos(d(w,u'))w}{\sin(d(w,u'))}$, we write $u' = \cos(d(w,u'))w + \sin(d(w,u'))(u'|\mathbb{S}_w)$ and, similarly, $u = \cos(d(w,u))w + \sin(d(w,u))(u|\mathbb{S}_w)$. Since $\mathbb{S}_w^\circ = \{\pm w\}$, we have $u'|\mathbb{S}_w = \mathrm{conv}(u', \{\pm w\}) \cap \mathbb{S}_w = u|\mathbb{S}_w$. Therefore,

$$\frac{u \cdot N_w^K}{\sin(d(u,w))} = (u|\mathbb{S}_w) \cdot N_w^K = \frac{u' \cdot N_w^K}{\sin(d(u',w))}.$$

Using Lemma 5.5 (iii) we conclude (6.6).

Now we can prove (6.4) for regular boundary points with positive curvature.
**Claim:** Let $w \in \operatorname{bd} K$ such that $H_{n-1}^{\mathbb{S}^n}(K, w) > 0$. Then
$$\lim_{\delta \to 0^+} \frac{-u \cdot N_w^K}{\sin(d(u,w))} \frac{d(w, w_\delta)}{\delta^{\frac{2}{n+1}}} = c_n H_{n-1}^{\mathbb{S}^n}(K, w)^{\frac{1}{n+1}}.$$

By the previous claim, we may assume that $K \subseteq C_w(\varepsilon)$ for arbitrarily small $\varepsilon > 0$ and $u \in \operatorname{int} K$ such that $K \subseteq \operatorname{int} \mathbb{S}_u^+$. Set $L = g_w(K)$, and $\xi = -\frac{g_w(u)}{\|g_w(u)\|}$. Since $\lim_{\delta \to 0^+} g_w(w_\delta) = 0$, we obtain
$$\lim_{\delta \to 0^+} \frac{-u \cdot N_w^K}{\sin(d(u,w))} \frac{d(w, w_\delta)}{\delta^{\frac{2}{n+1}}} = \lim_{\delta \to 0^+} \left(\xi \cdot N_w^K\right) \frac{\arctan(\|g_w(w_\delta)\|)}{\delta^{\frac{2}{n+1}}}$$
$$= \lim_{\delta \to 0^+} \left(\xi \cdot N_w^K\right) \frac{\|g_w(w_\delta)\|}{\delta^{\frac{2}{n+1}}}.$$

Using $\alpha = 0$, $\beta = \varepsilon$ and $u = w$ in Theorem 5.2, we conclude
$$\left\| x_{\frac{\delta}{\cos(\varepsilon)^{n+1}}} \right\| \geq \|g_w(w_\delta)\| \geq \|x_\delta\|$$
for $\delta$ small (note that the origin $0 = g_w(w)$ is a boundary point of $L$). Since $N_w^K = N_0^L$, Theorem 3.6 implies
$$\lim_{\delta \to 0^+} \left(\xi \cdot N_w^K\right) \frac{\|x_\delta\|}{\delta^{\frac{2}{n+1}}} = c_n H_{n-1}^{\mathbb{R}_w^n}(L, 0)^{\frac{1}{n+1}}$$
and, by Lemma 4.4, we have $H_{n-1}^{\mathbb{R}_w^n}(L, 0) = H_{n-1}^{\mathbb{S}^n}(K, w)$. We conclude
$$\frac{c_n}{\cos(\varepsilon)^2} H_{n-1}^{\mathbb{S}^n}(K, w)^{\frac{1}{n+1}} \geq \lim_{\delta \to 0^+} \frac{-u \cdot N_w^K}{\sin(d(u,w))} \frac{d(w, w_\delta)}{\delta^{\frac{2}{n+1}}} \geq c_n H_{n-1}^{\mathbb{S}^n}(K, w)^{\frac{1}{n+1}}.$$
Since $\varepsilon > 0$ can be chosen arbitrarily small the claim follows.

To finish the proof we only need to consider regular boundary points with vanishing curvature.
**Claim:** Let $K \subseteq C_u(\beta)$ for some $\beta \in \left(0, \frac{\pi}{2}\right)$ and $u \in \operatorname{int} K$. Then, for $w \in \operatorname{bd} K$ such that $H_{n-1}^{\mathbb{S}^n}(K, w) = 0$,
$$\lim_{\delta \to 0^+} \frac{-u \cdot N_w^K}{\sin(d(u,w))} \frac{d(w, w_\delta)}{\delta^{\frac{2}{n+1}}} = 0.$$

We consider $L = g_u(K)$ and $x = g_u(w)$. Then
$$\lim_{\delta \to 0^+} \frac{-u \cdot N_w^K}{\sin(d(u,w))} \frac{d(w, w_\delta)}{\delta^{\frac{2}{n+1}}} = \lim_{\delta \to 0^+} \left(\frac{x}{\|x\|} \cdot N_x^L\right) \frac{\arctan(\|x\|) - \arctan(\|g_u(w_\delta)\|)}{\delta^{\frac{2}{n+1}}}$$
$$\leq \lim_{\delta \to 0^+} \left(\frac{x}{\|x\|} \cdot N_x^L\right) \frac{\|x - g_u(w_\delta)\|}{\delta^{\frac{2}{n+1}}}.$$

By Theorem 5.2, we deduce
$$\|x - g_u(w_\delta)\| \leq \left\| x - x_{\frac{\delta}{\cos(\beta)^{n+1}}} \right\|.$$
As before, with Theorem 3.6 and Lemma 4.4, we conclude
$$0 \leq \lim_{\delta \to 0^+} \frac{-u \cdot N_w^K}{\sin(d(u,w))} \frac{d(w, w_\delta)}{\delta^{\frac{2}{n+1}}} \leq \lim_{\delta \to 0^+} \left(\frac{x}{\|x\|} \cdot N_x^L\right) \frac{\left\| x - x_{\frac{\delta}{\cos(\beta)^{n+1}}} \right\|}{\delta^{\frac{2}{n+1}}} = 0.$$
q.e.d.

This concludes the proof of Theorem 2.1.

## 7. The Floating Area

In this final section we will investigate some of the properties of the floating area. First we note that we may localize the floating area to a measure on the Borel $\sigma$-algebra $\mathcal{B}(\mathbb{S}^n)$ on $\mathbb{S}^n$ in the following way.

**Definition** (The Floating Measure)**.** *For $K \in \mathcal{K}(\mathbb{S}^n)$ and $\omega \in \mathcal{B}(\mathbb{S}^n)$ we define the* floating measure $\Omega(K,\omega)$ *by*

$$\Omega(K,\omega) = \int_{\omega \cap \operatorname{bd} K} H_{n-1}^{\mathbb{S}^n}(K,w)^{\frac{1}{n+1}}\, dw.$$

*The* floating area $\Omega(K)$ *of $K$ is given by* $\Omega(K) = \Omega(K, \mathbb{S}^n)$.

This notion is well defined since, by Theorem 2.1, the floating measure exists for all proper spherical convex bodies and is finite. For non-proper convex bodies the floating measure is identically zero. This is shown next.

**Theorem 7.1.** *The floating measure of a spherical polytope or a non-proper spherical convex body $C$ is trivial, that is, $\Omega(C, \cdot) \equiv 0$.*

*Proof.* This is obvious since in both cases the Gauss-Kronecker curvature of $C$ is zero for $\mathcal{H}^{n-1}$-almost all boundary points. Note that for $K = \mathbb{S}^n$ the floating measure is also trivial since $\operatorname{bd} \mathbb{S}^n = \emptyset$. q.e.d.

We first show that the floating measure is a valuation.

**Theorem 7.2.** *Let $K, L \in \mathcal{K}(\mathbb{S}^n)$ such that $K \cup L \in \mathcal{K}(\mathbb{S}^n)$. Then, for all $\omega \in \mathcal{B}(\mathbb{S}^n)$, we have*

$$\Omega(K,\omega) + \Omega(L,\omega) = \Omega(K \cup L, \omega) + \Omega(K \cap L, \omega).$$

*Proof.* This can be proofed similar to the Euclidean case by C. Schütt [**69**]. We will give a short outline of the argument. We decompose

$$\operatorname{bd}(K \cup L) = (\operatorname{bd} K \cap \operatorname{bd} L) \cup (\operatorname{bd} K \cap L^c) \cup (K^c \cap \operatorname{bd} L),$$
$$\operatorname{bd}(K \cap L) = (\operatorname{bd} K \cap \operatorname{bd} L) \cup (\operatorname{bd} K \cap \operatorname{int} L) \cup (\operatorname{int} L \cap \operatorname{bd} L),$$
$$\operatorname{bd} K = (\operatorname{bd} K \cap \operatorname{bd} L) \cup (\operatorname{bd} K \cap \operatorname{int} L) \cup (\operatorname{bd} K \cap L^c),$$
$$\operatorname{bd} L = (\operatorname{bd} K \cap \operatorname{bd} L) \cup (\operatorname{int} K \cap \operatorname{bd} L) \cup (K^c \cap \operatorname{bd} L),$$

where $K^c = \mathbb{S}^n \setminus K$ and $L^c = \mathbb{S}^n \setminus L$. Thus the integrals cancel for all sets but $\operatorname{bd} K \cap \operatorname{bd} L$. So we are done, once we show that

$$\int_{\operatorname{bd} K \cap \operatorname{bd} L} H_{n-1}^{\mathbb{S}^n}(K \cup L, w)^{\frac{1}{n+1}} dw + \int_{\operatorname{bd} K \cap \operatorname{bd} L} H_{n-1}^{\mathbb{S}^n}(K \cap L, w)^{\frac{1}{n+1}} dw$$
$$= \int_{\operatorname{bd} K \cap \operatorname{bd} L} H_{n-1}^{\mathbb{S}^n}(K, w)^{\frac{1}{n+1}} dw + \int_{\operatorname{bd} K \cap \operatorname{bd} L} H_{n-1}^{\mathbb{S}^n}(L, w)^{\frac{1}{n+1}} dw.$$

This follows from the fact that for $\mathcal{H}^{n-1}$-almost all $w \in \operatorname{bd} K \cap \operatorname{bd} L$ we have

$$H_{n-1}^{\mathbb{S}^n}(K \cup L, w) = \min\{H_{n-1}^{\mathbb{S}^n}(K,w), H_{n-1}^{\mathbb{S}^n}(L,w)\},$$
$$H_{n-1}^{\mathbb{S}^n}(K \cap L, w) = \max\{H_{n-1}^{\mathbb{S}^n}(K,w), H_{n-1}^{\mathbb{S}^n}(L,w)\}.$$

This local result follows from the Euclidean case by applying the gnomonic projection in $w$. q.e.d.

Next, we will show that the floating area is upper semicontinuous with respect to the spherical Hausdorff metric. The *spherical Hausdorff distance* $d^H$ of two convex bodies $K, L$ is defined by

$$d^H(K, L) = \min\{\lambda \geq 0 : K \subseteq L_\lambda \text{ and } L \subseteq K_\lambda\},$$

where $K_\lambda$ denotes the set of all points with distance at most $\lambda$ of $K$.

It is easy to see, that the floating area $\Omega(\cdot)$ cannot be continuous. We may consider a sequence of spherical polytopes $(P_j)_{j \in \mathbb{N}}$ that converges to a spherical cap $C_u(\frac{\pi}{4})$. Then $\Omega(C_u(\frac{\pi}{4})) = \omega_{n-1} 2^{-\frac{n-1}{2}} > 0$, but for all $j \in \mathbb{N}$ we have $\Omega(P_j, \mathbb{S}^n) = 0$.

**Theorem 7.3.** *The floating area is upper semicontinuous. Thus, for any sequence $(K_j)_{j \in \mathbb{N}}$ of convex bodies converging to a convex body $K$ in the spherical Hausdorff distance, we have $\limsup_{j \to \infty} \Omega(K_j, \mathbb{S}^n) \leq \Omega(K, \mathbb{S}^n)$.*

*Proof.* The proof of this theorem is along the lines of the proof of the upper semicontinuity of curvature integrals in the Euclidean setting by M. Ludwig [**39**]. Again, we will only give an outline. We denote the *m-th support measure* of a spherical convex body $K$ by $\Theta_m(K, \cdot)$. It is a uniquely determined finite Borel measure on $\mathbb{S}^n \times \mathbb{S}^n$ (see e.g. [**66**] for precise definitions). Furthermore it is weakly continuous in the first argument, that is, $K_j \to K$ in the Hausdorff metric implies $\Theta_m(K_j, \cdot) \xrightarrow{w} \Theta_m(K, \cdot)$.

The *m-th curvature measure* $C_m(K, \cdot)$ of $K \in \mathcal{K}(\mathbb{S}^n)$ is defined, for $\omega \in \mathcal{B}(\mathbb{S}^n)$, by $C_m(K, \omega) = \Theta_m(K, \omega \times \mathbb{S}^n)$. $C_m(K, \omega)$ is concentrated on $\operatorname{bd} K$. The Hausdorff measure restricted to $\operatorname{bd} K$ is denoted by $\mathcal{H}^{n-1}_{\operatorname{bd} K}$ and is a Radon measure. Thus we may split $C_m(K, \cdot) = C_m^a(K, \cdot) + C_m^s(K, \cdot)$ such that $C_m^a(K, \cdot)$ is absolute continuous with respect to $\mathcal{H}^{n-1}_{\operatorname{bd} K}$ and $C_m^s(K, \cdot)$ is the singular part. Moreover, for the absolute continuous part we have

$$C_m^a(K, \omega) = \int_{\omega \cap \operatorname{bd} K} H_{n-1-m}(K, w)\, dw,$$

for all $\omega \in \mathcal{B}(\mathbb{S}^n)$, and the singular part is concentrated on a null set, that is, there exists $\omega_0 \in \mathcal{B}(\mathbb{S}^n)$ such that $C_m^s(K, \omega \setminus \omega_0) = 0$ for all $\omega \in \mathcal{B}(\mathbb{S}^n)$.

Since $\Theta_m$ is weakly continuous in the first argument, so is $C_m$ and we conclude for $m = 0$,

$$\limsup_{j \to \infty} \int_{\omega \cap \operatorname{bd} K_j} H^{\mathbb{S}^n}_{n-1}(K_j, w)\, dw \leq \limsup_{j \to \infty} C_0(K_j, \omega) \leq C_0(K, \omega),$$

for all $\omega \in \mathcal{B}(\mathbb{S}^n)$. Using the arguments from [**39**, Section 4] and adapting the terminology, upper semicontinuity of the floating area follows. q.e.d.

Using the Gauss map $N^K_\cdot$, we find an equivalent expression for the floating area of $K$ as a curvature integral over the boundary of the polar of $K$.

**Theorem 7.4.** *Let $K \in \mathcal{K}(\mathbb{S}^n)$. Then*

$$\Omega(K) = \int_{\operatorname{bd} K^\circ} H^{\mathbb{S}^n}_{n-1}(K^\circ, w)^{\frac{n}{n+1}}\, dw.$$

*Proof.* The proof is similar to the proof of Theorem 2.8 in [**31**] by D. Hug. Even more can be said: Let $\omega \in \mathcal{B}(\mathbb{S}^n)$ and denote by $\sigma(K, \omega)$ the set of $v \in \mathbb{S}^n$ such that $v$ is an outer unit normal to some boundary point $w \in \operatorname{bd} K \cap \omega$. Clearly, $\sigma(K, \omega) \subseteq \operatorname{bd} K^\circ$. It follows that

$$\Omega(K, \omega) = \int_{\sigma(K,\omega)} H_{n-1}^{\mathbb{S}^n}(K^\circ, w)^{\frac{n}{n+1}} \, dw.$$

Since $\sigma(K, \mathbb{S}^n) = \operatorname{bd} K^\circ$, this implies the statement. q.e.d.

**7.1. Isoperimetric Inequality.** By Theorem 7.3, the floating area is upper semicontinuous on $\mathcal{K}(\mathbb{S}^n)$. Since $\mathcal{K}(\mathbb{S}^n)$ with the Hausdorff metric is compact (see, e.g., [**25**]), we may immediately conclude the existence of maximizers $C \in \mathcal{K}(\mathbb{S}^n)$ such that

$$\sup\{\Omega(K) : K \in \mathcal{K}(\mathbb{S}^n) \text{ such that } \operatorname{vol}_n(K) = c\} \leq \Omega(C)$$

for a fixed $c \in [0, \frac{\omega_n}{2}]$ and $\operatorname{vol}_n(C) = c$.

We believe that the only maximizers of the floating area are geodesic balls:

**Conjecture 7.5.** Let $K \in \mathcal{K}(\mathbb{S}^n)$. Then

$$\Omega(K) \leq \Omega(C^K), \tag{7.1}$$

where $C^K$ is a spherical cap such that $\operatorname{vol}_n(C^K) = \operatorname{vol}_n(K)$. Equality holds if and only if $K$ is spherical cap.

This conjecture is still open, but we are able to prove the following inequality.

**Theorem 7.6.** Let $K \in \mathcal{K}(\mathbb{S}^n)$. Then

$$\Omega(K) \leq P(K^\circ)^{\frac{1}{n+1}} P(K)^{\frac{n}{n+1}}, \tag{7.2}$$

where $P(K)$ denotes the perimeter of $K$. Equality holds if and only if $K$ is a spherical cap.

*Proof.* Using Hölder's inequality, we obtain

$$\Omega(K) = \int_{\operatorname{bd} K} H_{n-1}^{\mathbb{S}^n}(K, w)^{\frac{1}{n+1}} \, dw \leq \left( \int_{\operatorname{bd} K} H_{n-1}^{\mathbb{S}^n}(K, w) \, dw \right)^{\frac{1}{n+1}} \left( \int_{\operatorname{bd} K} dw \right)^{\frac{n}{n+1}}.$$

Since $N^K(\operatorname{bd} K) \subseteq \operatorname{bd} K^\circ$ and $J^{\operatorname{bd} K}(N^K)(w) = H_{n-1}^{\mathbb{S}^n}(K, w)$ for $\mathcal{H}^{n-1}$-almost all $w \in \operatorname{bd} K$, this implies (7.2).

That equality holds precisely for spherical caps follows from the fact, that equality holds in Hölder's inequality if and only if $H_{n-1}^{\mathbb{S}^n}(K, .)$ is constant. q.e.d.

Another inequality for the floating area can be derived using the gnomonic projection and the affine isoperimetric inequality.

**Theorem 7.7.** Let $K \in \mathcal{K}_0^p(\mathbb{S}^n)$. Then, for $u \in \operatorname{int} K$ and $0 < \alpha \leq \beta < \frac{\pi}{2}$ such that $C_u(\alpha) \subseteq K \subseteq C_u(\beta)$, we have

$$\frac{\Omega(K)}{\omega_{n-1}} \leq \left( \frac{\cos(\alpha)^2 \tan(\beta)^n}{\tan(\alpha)^{n-1}} \frac{P(K)}{\omega_{n-1}} \right)^{\frac{n-1}{n+1}}. \tag{7.3}$$

*Equality holds for spherical caps for which $\alpha = \beta$.*

*Proof.* Using the gnomonic projection in $u$, (4.11) and (4.13), we obtain

$$\Omega(K) = \int_{\operatorname{bd} g_u(K)} H_{n-1}^{\mathbb{S}^n}(K, g_u^{-1}(x))^{\frac{1}{n+1}} J^{\operatorname{bd} g_u(K)}(g_u^{-1})(x) \, dx$$

$$= \int_{\operatorname{bd} g_u(K)} H_{n-1}^{\mathbb{R}_u^n}(g_u(K), x)^{\frac{1}{n+1}} \frac{1}{(1 + \|x\|^2)^{\frac{n-1}{2}}} \, dx.$$

Since $\tan(\alpha) \leq \|x\| \leq \tan(\beta)$ for all $x \in \operatorname{bd} g_u(K)$, we conclude
$$\cos(\beta)^{n-1} \operatorname{as}(g_u(K)) \leq \Omega(K) \leq \cos(\alpha)^{n-1} \operatorname{as}(g_u(K)).$$

Using the classical affine isoperimetric inequality for convex bodies $L \in \mathcal{K}_0(\mathbb{R}^n)$,
$$\frac{\operatorname{as}(L)}{\omega_{n-1}} \leq \left(\frac{n \operatorname{vol}_n(L)}{\omega_{n-1}}\right)^{\frac{n-1}{n+1}},$$
gives
$$\frac{\Omega(K)}{\omega_{n-1}} \leq \left(\cos(\alpha)^{n+1} \frac{n \operatorname{vol}_n(g_u(K))}{\omega_{n-1}}\right)^{\frac{n-1}{n+1}}.$$

For the volume of the gnomonic projection we derive the inequality
$$\operatorname{vol}_n(g_u(K)) \stackrel{(4.7)}{=} \int_K \frac{dv}{\cos(d(v,u))^{n+1}} = \int_{\operatorname{bd} K} \frac{-u \cdot N_w^K}{\sin(d(u,w))^n} \int_0^{d(u,w)} \frac{\tan(t)^{n-1}}{\cos(t)^2} \, dt$$
$$\leq \int_{\operatorname{bd} K} \frac{1}{\sin(d(u,w))^{n-1}} \frac{\tan(d(u,w))^n}{n} \, dw \leq \frac{P(K)}{n} \frac{\tan(\beta)^n}{\sin(\alpha)^{n-1}},$$
where we used the fact that $-u \cdot N_w^K = \sin(d(u, \mathbb{S}_{N_w^K})) \leq \sin(d(u,w))$. This concludes the prove of (7.3). It is easy to check that equality holds for spherical caps of radius $\alpha = \beta$. q.e.d.

Inequality (7.3) is weaker than our conjectured inequality (7.1): For any $K \in \mathcal{K}_0^p(\mathbb{S}^n)$ and $u \in \operatorname{int} K$ such that $C_u(\alpha) \subseteq K \subseteq C_u(\beta)$ we have
$$\operatorname{vol}_n(C_u(\alpha)) \leq \operatorname{vol}_n(K) \leq \operatorname{vol}_n(C_u(\beta)).$$

Thus, for the spherical cap $C^K = C_u(\alpha_K)$ such that $\operatorname{vol}_n(C^K) = \operatorname{vol}_n(K)$, we have $\alpha \leq \alpha_K \leq \beta$. We therefore conclude for the right hand side of inequality (7.3) that
$$\left(\cos(\alpha_K) \sin(\alpha_K) \frac{P(K)}{\omega_{n-1}}\right)^{\frac{n-1}{n+1}} \leq \left(\frac{\cos(\alpha)^2 \tan(\beta)^n}{\tan(\alpha)^{n-1}} \frac{P(K)}{\omega_{n-1}}\right)^{\frac{n-1}{n+1}}. \qquad (7.4)$$

Our conjectured inequality (7.1) would imply
$$\Omega(K) \leq \Omega(C^K) = \cos(\alpha_K)^{\frac{n-1}{n+1}} \sin(\alpha_K)^{n \frac{n-1}{n+1}} \omega_{n-1},$$
which in turn would imply (7.3) by (7.4) and the isoperimetric inequality
$$P(K) \geq P(C^K) = \sin(\alpha_K)^{n-1} \omega_{n-1}.$$

**Acknowledgement.** The authors would like to thank the Institute for Mathematics and Applications (IMA), University of Minnesota. It was during their stay there that part of the paper was written

## References

[1] J. Abardia, E. Gallego, and G. Solanes, *The Gauss-Bonnet theorem and Crofton-type formulas in complex space forms*, Israel J. Math. **187** (2012), 287–315.

[2] S. Alesker, *Theory of valuations on manifolds: a survey*, Geom. Funct. Anal. **17** (2007), no. 4, 1321–1341.

[3] S. Alesker, *Valuations on manifolds and integral geometry*, Geom. Funct. Anal. **20** (2010), no. 5, 1073–1143.

[4] D. Amelunxen, M. Lotz, M. B. McCoy, and J. A. Tropp, *Living on the edge: A geometric theory of phase transitions in convex optimization*, CoRR (2013).

[5] B. Andrews, *Contraction of convex hypersurfaces by their affine normal*, J. Differential Geom. **43** (1996), no. 2, 207–230.

[6] B. Andrews, *The affine curve-lengthening flow*, J. Reine Angew. Math. **506** (1999), 43–83.

[7] S. Artstein-Avidan, B. Klartag, C. Schütt, and E. Werner, *Functional affine-isoperimetry and an inverse logarithmic Sobolev inequality*, J. Funct. Anal. **262** (2012), no. 9, 4181–4204.

[8] A. Bernig, J. H. G. Fu, and G. Solanes, *Integral geometry of complex space forms*, Geom. Funct. Anal. **24** (2014), no. 2, 403–492.

[9] F. Besau and F. E. Schuster, *Binary Operations in Spherical Convex Geometry*, ArXiv e-prints (July 2014), available at `1407.1153`.

[10] W. Blaschke, *Vorlesung über Differentialgeometrie II, Affine Differntialgeometrie*, Springer-Verlag, Berlin, 1923.

[11] K. Böröczky Jr., *Approximation of general smooth convex bodies*, Adv. Math. **153** (2000), no. 2, 325–341.

[12] K. Böröczky Jr., *Polytopal approximation bounding the number of k-faces*, J. Approx. Theory **102** (2000), no. 2, 263–285.

[13] K. J. Böröczky, *Stability of the Blaschke-Santaló and the affine isoperimetric inequality*, Adv. Math. **225** (2010), no. 4, 1914–1928.

[14] U. Caglar and E. M. Werner, *Divergence for s-concave and log concave functions*, Adv. Math. **257** (2014), 219–247.

[15] C. Dupin, *Application de géométrie et de méchanique*, Paris, 1822.

[16] H. Federer, *Geometric measure theory*, Die Grundlehren der mathematischen Wissenschaften, Band 153, Springer-Verlag New York Inc., New York, 1969.

[17] A. Figalli and Y. Ge, *Isoperimetric-type inequalities on constant curvature manifolds*, Adv. Calc. Var. **5** (2012), no. 3, 251–284.

[18] T. Figiel, J. Lindenstrauss, and V. D. Milman, *The dimension of almost spherical sections of convex bodies*, Acta Math. **139** (1977), no. 1-2, 53–94.

[19] B. Fleury, O. Guédon, and G. Paouris, *A stability result for mean width of $L_p$-centroid bodies*, Adv. Math. **214** (2007), no. 2, 865–877.

[20] F. Gao, D. Hug, and R. Schneider, *Intrinsic volumes and polar sets in spherical space*, Math. Notae **41** (2003), 159–176.

[21] R. J. Gardner and G. Zhang, *Affine inequalities and radial mean bodies*, Amer. J. Math. **120** (1998), no. 3, 505–528.

[22] R. J. Gardner, *Geometric tomography*, Second, Encyclopedia of Mathematics and its Applications, vol. 58, Cambridge University Press, Cambridge, 2006.

[23] C. Gerhardt, *Curvature flows in the sphere*, J. Differential Geom. (to appear), available at `1308.1607`.

[24] C. Gerhardt, *Minkowski type problems for convex hypersurfaces in the sphere*, Pure Appl. Math. Q. **3** (2007), no. 2, part 1, 417–449.

[25] S. Glasauer, *Integral geometry of spherically convex bodies*, Diss. Summ. Math. **1** (1996), no. 1-2, 219–226.

[26] P. M. Gruber, *Volume approximation of convex bodies by inscribed polytopes*, Math. Ann. **281** (1988), no. 2, 229–245.

[27] P. M. Gruber, *Aspects of approximation of convex bodies*, Handbook of convex geometry, Vol. A, B, 1993, pp. 319–345.

[28] P. M. Gruber, *Convex and discrete geometry*, Grundlehren der Mathematischen Wissenschaften [Fundamental Principles of Mathematical Sciences], vol. 336, Springer, Berlin, 2007.

[29] C. Haberl and L. Parapatits, *The centro-affine Hadwiger theorem*, J. Amer. Math. Soc. **27** (2014), no. 3, 685–705.

[30] C. Haberl and F. E. Schuster, *General $L_p$ affine isoperimetric inequalities*, J. Differential Geom. **83** (2009), no. 1, 1–26.

[31] D. Hug, *Contributions to affine surface area*, Manuscripta Math. **91** (1996), no. 3, 283–301.

[32] M. N. Ivaki, *On the stability of the p-affine isoperimetric inequality*, J. Geom. Anal. **24** (2014), no. 4, 1898–1911.

[33] M. N. Ivaki and A. Stancu, *Volume preserving centro-affine normal flows*, Comm. Anal. Geom. **21** (2013), no. 3, 671–685.

[34] D. A. Klain and G.-C. Rota, *Introduction to geometric probability*, Lezioni Lincee. [Lincei Lectures], Cambridge University Press, Cambridge, 1997.

[35] S. G. Krantz and H. R. Parks, *Geometric integration theory*, Cornerstones, Birkhäuser Boston, Inc., Boston, MA, 2008.

[36] K. Leichtweiß, *Zur Affinoberfläche konvexer Körper*, Manuscripta Math. **56** (1986), no. 4, 429–464.

[37] K. Leichtweiß, *Affine geometry of convex bodies*, Johann Ambrosius Barth Verlag, Heidelberg, 1998.

[38] M. Ludwig, *Asymptotic approximation of smooth convex bodies by general polytopes*, Mathematika **46** (1999), no. 1, 103–125.

[39] M. Ludwig, *On the semicontinuity of curvature integrals*, Math. Nachr. **227** (2001), 99–108.

[40] M. Ludwig, *General affine surface areas*, Adv. Math. **224** (2010), no. 6, 2346–2360.

[41] M. Ludwig and M. Reitzner, *A characterization of affine surface area*, Adv. Math. **147** (1999), no. 1, 138–172.

[42] M. Ludwig and M. Reitzner, *A classification of* SL($n$) *invariant valuations*, Ann. of Math. (2) **172** (2010), no. 2, 1219–1267.

[43] E. Lutwak, *On some affine isoperimetric inequalities*, J. Differential Geom. **23** (1986), no. 1, 1–13.

[44] E. Lutwak, *Extended affine surface area*, Adv. Math. **85** (1991), no. 1, 39–68.

[45] E. Lutwak, *The Brunn-Minkowski-Firey theory. II. Affine and geominimal surface areas*, Adv. Math. **118** (1996), no. 2, 244–294.

[46] E. Lutwak and V. Oliker, *On the regularity of solutions to a generalization of the Minkowski problem*, J. Differential Geom. **41** (1995), no. 1, 227–246.

[47] E. Lutwak, D. Yang, and G. Zhang, *$L_p$ affine isoperimetric inequalities*, J. Differential Geom. **56** (2000), no. 1, 111–132.

[48] E. Lutwak, D. Yang, and G. Zhang, *The Cramer-Rao inequality for star bodies*, Duke Math. J. **112** (2002), no. 1, 59–81.

[49] E. Lutwak, D. Yang, and G. Zhang, *Moment-entropy inequalities*, Ann. Probab. **32** (2004), no. 1B, 757–774.

[50] E. Lutwak, D. Yang, and G. Zhang, *Cramér-Rao and moment-entropy inequalities for Renyi entropy and generalized Fisher information*, IEEE Trans. Inform. Theory **51** (2005), no. 2, 473–478.

[51] E. Lutwak and G. Zhang, *Blaschke-Santaló inequalities*, J. Differential Geom. **47** (1997), no. 1, 1–16.

[52] F. Maggi, *Sets of finite perimeter and geometric variational problems*, Cambridge Studies in Advanced Mathematics, vol. 135, Cambridge University Press, Cambridge, 2012. An introduction to geometric measure theory.

[53] M. Makowski and J. Scheuer, *Rigidity results, inverse curvature flows and Alexandrov-Fenchel type inequalities in the sphere*, ArXiv e-prints (July 2013), available at 1307.5764.

[54] M. Meyer and E. Werner, *On the p-affine surface area*, Adv. Math. **152** (2000), no. 2, 288–313.

[55] K. Nomizu and T. Sasaki, *Affine differential geometry*, Cambridge Tracts in Mathematics, vol. 111, Cambridge University Press, Cambridge, 1994. Geometry of affine immersions.

[56] V. Oliker, *Embedding* $\mathbf{S}^n$ *into* $\mathbf{R}^{n+1}$ *with given integral Gauss curvature and optimal mass transport on* $\mathbf{S}^n$, Adv. Math. **213** (2007), no. 2, 600–620.

[57] G. Paouris and E. M. Werner, *Relative entropy of cone measures and $L_p$ centroid bodies*, Proc. Lond. Math. Soc. (3) **104** (2012), no. 2, 253–286.

[58] C. M. Petty, *Affine isoperimetric problems*, Discrete geometry and convexity (New York, 1982), 1985, pp. 113–127.

[59] M. Reitzner, *The combinatorial structure of random polytopes*, Adv. Math. **191** (2005), no. 1, 178–208.

[60] L. A. Santaló, *An affine invariant for convex bodies of n-dimensional space*, Portugaliae Math. **8** (1949), 155–161.

[61] L. A. Santaló, *On parallel hypersurfaces in the elliptic and hyperbolic n-dimensional space*, Proc. Amer. Math. Soc. **1** (1950), 325–330.

[62] L. A. Santaló, *Cauchy and Kubota's formula for convex bodies in elliptic n-space*, Rend. Sem. Mat. Univ. Politec. Torino **38** (1980), no. 1, 51–58.

[63] L. A. Santaló, *Integral geometry and geometric probability*, Second, Cambridge Mathematical Library, Cambridge University Press, Cambridge, 2004. With a foreword by Mark Kac.

[64] R. Schneider, *The endomorphisms of the lattice of closed convex cones*, Beiträge Algebra Geom. **49** (2008), no. 2, 541–547.

[65] R. Schneider, *Convex bodies: the Brunn-Minkowski theory*, expanded, Encyclopedia of Mathematics and its Applications, vol. 151, Cambridge University Press, Cambridge, 2014.

[66] R. Schneider and W. Weil, *Stochastic and integral geometry*, Probability and its Applications (New York), Springer-Verlag, Berlin, 2008.

[67] F. E. Schuster and T. Wannerer, GL($n$) *contravariant Minkowski valuations*, Trans. Amer. Math. Soc. **364** (2012), no. 2, 815–826.

[68] C. Schütt, *The convex floating body and polyhedral approximation*, Israel J. Math. **73** (1991), no. 1, 65–77.

[69] C. Schütt, *On the affine surface area*, Proc. Amer. Math. Soc. **118** (1993), no. 4, 1213–1218.

[70] C. Schütt, *Random polytopes and affine surface area*, Math. Nachr. **170** (1994), 227–249.

[71] C. Schütt and E. Werner, *The convex floating body*, Math. Scand. **66** (1990), no. 2, 275–290.

[72] C. Schütt and E. Werner, *Polytopes with vertices chosen randomly from the boundary of a convex body*, Geometric aspects of functional analysis, 2003, pp. 241–422.

[73] C. Schütt and E. Werner, *Surface bodies and p-affine surface area*, Adv. Math. **187** (2004), no. 1, 98–145.

[74] G. Solanes, *Integral geometry and the Gauss-Bonnet theorem in constant curvature spaces*, Trans. Amer. Math. Soc. **358** (2006), no. 3, 1105–1115 (electronic).

[75] A. Stancu, *The discrete planar $L_0$-Minkowski problem*, Adv. Math. **167** (2002), no. 1, 160–174.

[76] A. Stancu, *On the number of solutions to the discrete two-dimensional $L_0$-Minkowski problem*, Adv. Math. **180** (2003), no. 1, 290–323.

[77] N. S. Trudinger and X.-J. Wang, *The affine Plateau problem*, J. Amer. Math. Soc. **18** (2005), no. 2, 253–289.

[78] E. Werner and D. Ye, *New $L_p$ affine isoperimetric inequalities*, Adv. Math. **218** (2008), no. 3, 762–780.

[79] E. M. Werner, *Rényi divergence and $L_p$-affine surface area for convex bodies*, Adv. Math. **230** (2012), no. 3, 1040–1059.

[80] V. Yaskin, *The Busemann-Petty problem in hyperbolic and spherical spaces*, Adv. Math. **203** (2006), no. 2, 537–553.

*E-mail address*: florian.besau@tuwien.ac.at

(F. Besau)

INSTITUTE OF DISCRETE MATHEMATICS AND GEOMETRY
VIENNA UNIVERSITY OF TECHNOLOGY
WIEDNER HAUPTSTRASSE 8–10/1046
1040 WIEN, AUSTRIA

*E-mail address*: elisabeth.werner@case.edu

(E. M. Werner)

DEPARTMENT OF MATHEMATICS
CASE WESTERN RESERVE UNIVERSITY
CLEVELAND, OHIO 44106, U. S. A.